\documentclass[pdflatex,sn-vancouver]{sn-jnl}

\jyear{2023}

\usepackage{subfigure}
\usepackage{paracol}
\usepackage{natbib}
\usepackage[font=small,labelfont=bf,tableposition=top]{caption}
\usepackage{graphics}
\usepackage{afterpage} 

\theoremstyle{thmstyleone}%
\usepackage{afterpage}
\usepackage{multicol,multirow,tabulary}
\usepackage{float} 

\usepackage[section]{placeins} 

\begin{document}

\title[Congressional Districting: ``Rocks-Pebbles-Sand'' Approach]{A Case Study of Congressional Districting:\\ ``Rocks-Pebbles-Sand'' Approach}


\author[1]{\sur{\href{https://scholar.google.com/citations?user=L94aiQ4AAAAJ&hl=en}{Jimmy Risk}}}\email{jrisk@cpp.edu}

\author[1]{\sur{\href{https://scholar.google.com/citations?user=YfO2FToAAAAJ&hl=en&oi=ao}{Jennifer Switkes}}}\email{jmswitkes@cpp.edu}

\author[1]{\fnm{Ann} \sur{Zhang}}\email{yuz2@cpp.edu}

\affil[1]{\orgdiv{Department of Mathematics and Statistics}, \orgname{California State Polytechnic University, Pomona}, \orgaddress{\street{3801 W. Temple Avenue}, \city{Pomona},  \state{California}, \postcode{91768},\country{U.S.A}}}

\abstract{As a case study into an algorithmic approach to congressional districting, North Carolina provides a lot to explore. Statistical modeling has called into question whether recent North Carolina district plans are unbiased. In particular, the literature suggests that the district plan used in the 2016 U.S. House of Representatives election yields outlier results that are statistically unlikely to be obtained without the application of bias. Therefore, methods for creating strong and fair district plans are needed. Informed by previous districting models and algorithms, we build a model and algorithm to produce an ensemble of viable Congressional district plans.  
Our work contributes a ``Rocks-Pebbles-Sand'' concept and procedure facilitating reasonable population equity with a small overall number of county splits among districts. Additionally, our methodology minimizes the initial need for granular, precinct-level data, thereby reducing the risk of covert gerrymandering. This case study indicates plausibility of an approach built upon an easy-to-understand intuition.
}

\keywords{Redistricting, case study, algorithm, district plan, gerrymandering}


\maketitle

\section{Introduction}

In the complex landscape of congressional redistricting, North Carolina serves as a complex and instructive case study. Our work aligns with existing algorithms in its foundational structure, notably those discussed by \citet{cirincione2000assessing}, which also build districts iteratively from smaller units. However, our algorithm departs from uniformity by segregating its operation into discrete phases. Through a hierarchical ``rocks-pebbles-sand" strategy, we offer a structured yet nuanced approach to district formation.

A noteworthy advantage of our methodology is its minimization of dependency on granular, precinct-level data, thereby reducing the risk of insidious gerrymandering. Although precinct-based mapping is traditional and widespread, our approach demonstrates that effective districting can commence with fewer, larger building blocks. This facilitates a simplified initial map construction with the option for further granularity in later stages, all while maintaining algorithmic integrity.

Designed for a general audience, our work not only presents a distinct methodology but also underscores the importance of merging pre-split counties, aligning with the broader aim of creating fair and unbiased district plans. This fairness is particularly critical in a state like North Carolina, where recent election outcomes have raised questions about the balance of political representation, as discussed below.

From $2014$ to $2020$, North Carolina was divided into $13$ Congressional districts. In $2014$, $2016$, and $2018$, Democrats won $3$ out of $13$ seats and Republicans won $10$ out of $13$ seats. In $2020$, Democrats won $5$ out of $13$ seats and Republicans won $8$ out of $13$ seats. In Table~\ref{tab: Results of the United States House of Representatives Elections in North Carolina from 2016 - 2020}, we show the results of recent United States House of Representatives elections in North Carolina. Since the parties have close-to-balanced vote outcomes in each election, one might wonder if the unbalanced result is due to political gerrymandering or is a natural consequence of North Carolina's spatial partisan distribution and the associated geopolitical structure. 

\begin{table}[ht]
    \centering
    \begin{tabular}{|c|c|c|c|c|}
    \hline
    Election Year     & Seats won  & Seats won  & Democrat & Republican\\  
& by Democrat & by Republican & \multicolumn{1}{c|}{Votes (\%)} & \multicolumn{1}{c|}{Votes (\%)} \\ \hline
      2014 & 3 & 10 & 44.0\% & 55.4\% \\\hline
      2016 & 3 & 10 & 46.6\% & 53.2\% \\\hline
      2018 & 3 & 10 & 48.3\% & 50.5\%\\\hline
      2020 & 5 & 8  & 50.0\% & 49.4\% \\\hline
    \end{tabular}
    \caption{Results of the United States House of Representatives Elections in North Carolina from 2014 - 2020 (Source: State of North Carolina Election Data, \url{https://er.ncsbe.gov/}).}
    \label{tab: Results of the United States House of Representatives Elections in North Carolina from 2016 - 2020}
\end{table}

Every ten years, states redraw Congressional district lines based on census data. The redistricting of each state is determined by the following criteria: population equity, racial and language minority protections, compactness and contiguity, political subdivisions, and communities of interest (\citet{congressional}). The counties or sub-counties of each district should be geographically contiguous, meaning that one can walk from any point to any other point in the district without going outside of the district.
Under the ``One Person--One Vote Principle,'' the population size of each district  should be equal or nearly equal, especially in a majoritarian system (\citet{apollonio1}).
Each district should be closely and neatly packed together. Various measures of district compactness appear in the literature.

\begin{figure}
    \centering
    \includegraphics[scale=0.6]{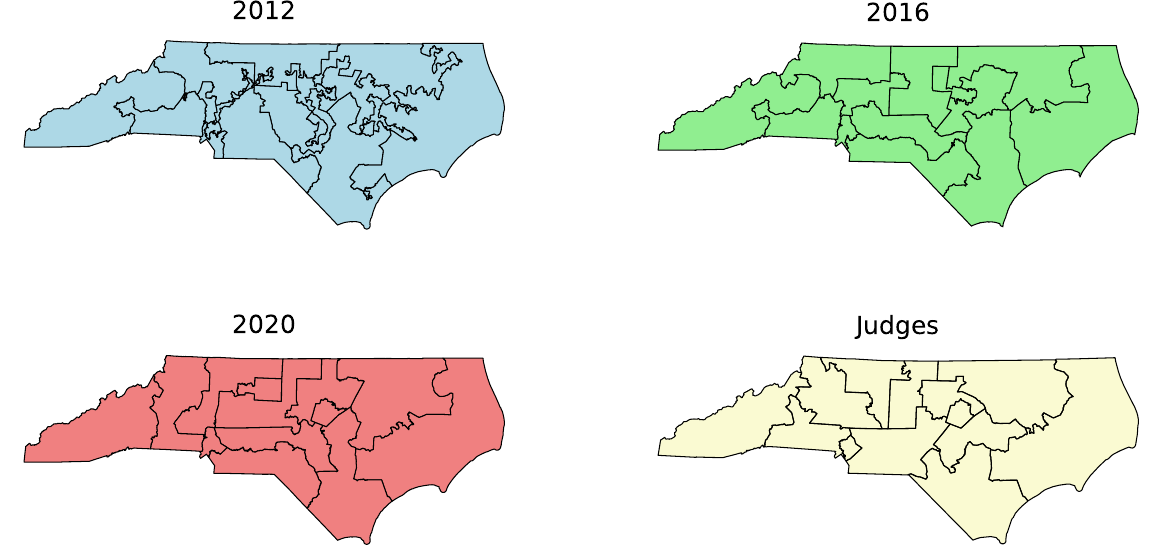}
    \caption{Congressional District Plan of North Carolina (Source: North Carolina General Assembly, \url{https://www.ncleg.gov/Redistricting}).}
    \label{fig:historic_plans}
\end{figure}

The North Carolina Congressional district map that was used for the $2012$, $2016$, and $2020$ elections is shown in Figure~\ref{fig:historic_plans}. Although the newer Congressional map no longer had a snake-like shaped district like was present in $2014$, the outcome in terms of seats won in $2016$ and $2018$ did not change after the new Congressional map was introduced. Also shown in Figure~\ref{fig:historic_plans} is a simulated district map created by a panel of retired judges in 2016 as part of a research project into fair redistricting; in Section~\ref{sec:analysis}, this map is compared with our simulations and the official district maps.

We dive a little deeper into the 2016 election. 

\begin{itemize}
    \item In 2016, Democrats won the seats in Districts 1, 4, and 12 with close to $70\%$ of the votes in each corresponding district. The margin of victory for Democrats was at least $34\%$ in each of these three districts. On the other hand, the margin of victory for six Republican districts was under $20\%$. By packing the votes of Democrats into certain areas, map makers may have reduced the influence of votes by Democrats in other districts.
    
   \item Wake county was split into two sub-counties, with one belonging to District 2 and the other one belonging to District 4. The Republicans won District 2, and the part of Wake county that belongs to District 2 happens to capture all the Republican-leaning precincts within its uncommon shape. Meanwhile, the Democrats won District 4 with about a $36\%$ margin of victory. If we were to split the original Wake county in two by using a roughly horizontal or vertical line, then the Democrats would have a much higher chance of winning both districts instead of one.
    
   \item Mecklenburg county was split into two sub-counties. The first sub-county belonged to District 9 and it consisted of about $23\%$ of the votes from Mecklenburg county. The second sub-county belonged to District 12 and it consisted of  about $77\%$ of the votes. In the first sub-county, Democrats had fewer votes than Republicans, but votes for Democrats were more than twice as many as votes for Republicans in the second sub-county.

\end{itemize}

In summary, in 2016 the districts that Democrats won tended to have a very large margin of victory, an indicator of possible districting bias.

However, assessing gerrymandering is not an easy task as the redistricting criteria are difficult to rigorously define and assess.  Furthermore, the political nature of the redistricting process naturally causes disagreement.  Looking at the historical evolution of assessment of gerrymandering, the Supreme Court first found partisan gerrymandering to be justiciable in the \emph{Bandemer v.~Davis (1985)} case (\citet{geldzahler1988davis}).  In \emph{Veith v.~Jubelirer (2004)}, five out of nine justices agreed that some standard of gerrymandering might be adopted in the future if such a standard can be found (\citet{grofman2007future}).  \citet{godfrey2006brief} proposed \emph{partisan symmetry} which principles two parties to receive the same percentage of seats if they obtain the same percentage of votes.  Following the work of \citet{godfrey2006brief}, there was strong agreement in the scholarly community that partisan symmetry was the de facto method to assess gerrymandering (\citet{stephanopoulos2018measure}).  However, in a 2006 case, Justice Kennedy expressed concern that such a measure ``may in large depend on conjecture about where vote-switchers will reside" (\citet{stephanopoulos2018measure}).  Indeed, while partisan symmetry is one useful measure, it has its flaws (see \citet{deford2021implementing} for a thorough investigation).  The interested reader can find further discussion and alternative mathematical approaches to gerrymandering in, e.g., \citet{wang2016three} and 
\citet{stephanopoulos2018measure}.  Thus, there is a difficulty in a consensus and evaluation of a general gerrymandering criteria. 

On the state level, North Carolina passed House Bill 92 (HB92) in 2015.  This bill is a step in the right direction to establish a nonpartisan redistricting process.  It provides the following redistricting standards: (i) population equity: ``districts shall each have a population as nearly equal as practicable to the ideal population," (ii) minimizing county splits: ``the number of counties and cities divided among more than one district shall be as small as possible," (iii) contiguity, (iv) compactness, referencing \textit{length-width compactness} and \textit{perimeter compactness} measures, and (v) concordance with the Voting Rights Act of 1965 (VRA).  

In addition, specific partisan metrics have been used in North Carolina court cases, such as the \emph{efficiency gap} and \emph{mean-median difference measure} used in \textit{Common Cause v. Rucho} (2018) to determine the unconstitutional partisan gerrymander of the 2016 congressional district plan. \textit{Common Cause v. Lewis} (2019) utilized the \emph{lopsided margins test} in determining the violation of a state legislative district map. These partisan measures illustrate that partisan bias can exist despite passing traditional metrics. In particular, the 2016 plan was generally adequate in terms of compactness, which is one reason why its unconstitutional partisan gerrymander was not identified until after the fact.

Our main goal is to develop a redistricting algorithm based on intuitive notions that appeal to a broader audience. 

\begin{figure}[h!]
    \centering
    \includegraphics[scale=0.45]{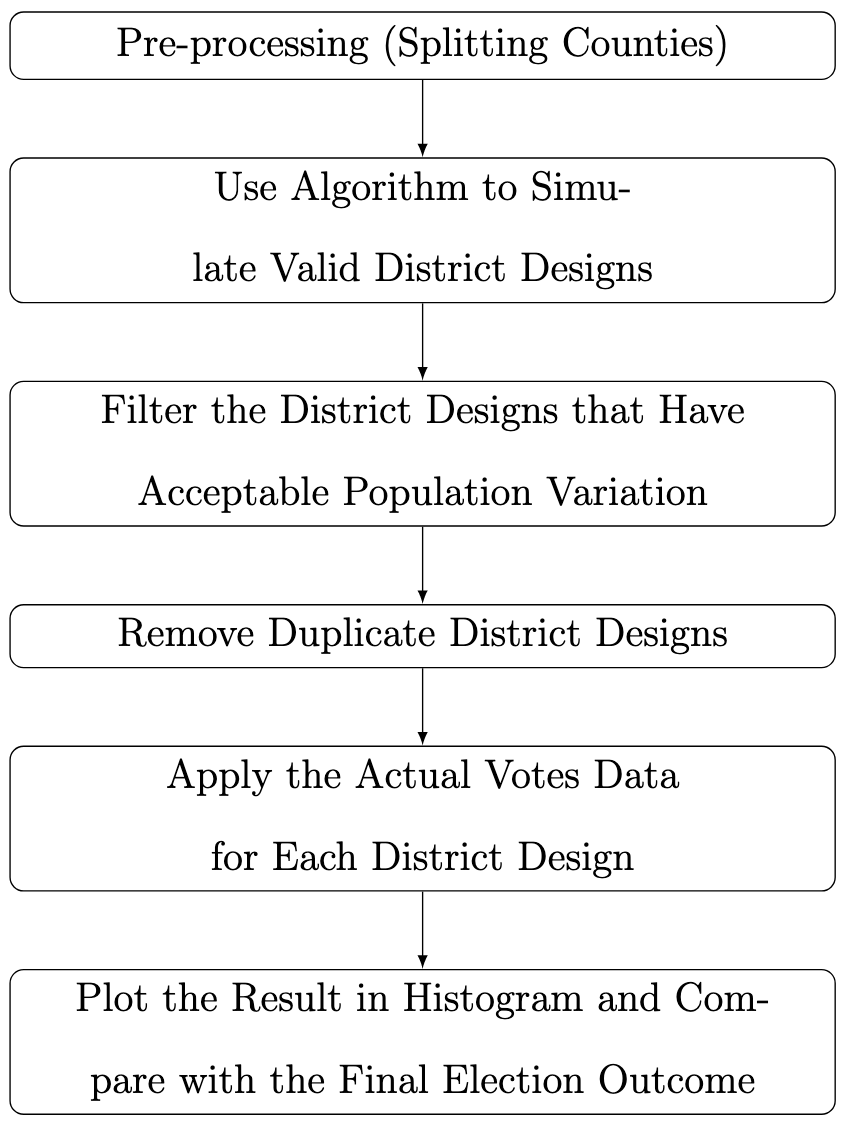}
    \caption{Flowchart for North Carolina District Plan Model.}
    \label{flowchart}
\end{figure}

While modern computing allows for the creation of sophisticated plans using mathematical measures, we find that only a small number of county splits are necessary to produce plausible results based on the literature and 2015 North Carolina House Bill 92 (HB92). This task of minimizing county splits is also explored by Shahmizad et al. \citep{shahmizad2023political} through integer programming. Although metrics like \emph{split pairs} suggested by \citet{waschpress2021split} exist for evaluating county splits, they did not offer additional insights in our context due to the small number of county splits in our plans. Our algorithm enables us to construct redistricting plans that meet provided constraints while being straightforward and transparent. This is particularly relevant as the mathematics of gerrymandering can be convoluted and enigmatic. An overview of our approach is presented in Figure~\ref{flowchart}.  

To elucidate on existing redistricting algorithms, there are three main classifications according to \citet{herschlag2020quantifying}. \emph{Constructive randomized algorithms} generate redistricting plans with an initial random seed and either grow a fixed number of districts or combine small districts until the number of desired districts is achieved (e.g., \citet{cirincione2000assessing, chen2013unintentional, chen2015cutting}). \emph{Moving boundary MCMC algorithms} begin with a redistricting plan and modify the boundaries according to a target distribution encompassing plausible redistricting criteria (e.g., \citet{macmillan2001redistricting,  mattingly2014redistricting, bangia2015quantifying, wu2015impartial, imai2015new, herschlag2020quantifying}). \emph{Optimization algorithms} seek optimal redistricting plans on a highly non-convex space (e.g., \citet{mehrotra1998optimization, liu2016pear}).  

While previous discussions have emphasized the need for a more mathematical approach to gerrymandering, there is also a clear reluctance to fully adopt complex methods. Our method, while bearing similarities to constructive randomized algorithms like those of \citet{cirincione2000assessing}, introduces a hierarchical strategy that allows for a nuanced yet structured formation of districts without node removal. Our approach is based on an intuitive notion, using a ``Rocks-Pebbles-Sand'' concept that prioritizes population equity while minimizing county splits. This approach can pass common gerrymandering tests and guarantees fewer county splits compared to existing plans. Additionally, we find that introducing constraints can help address population equity, though it may be difficult to satisfy the exact 0.1\% constraint mandated by HB92. For a detailed discussion on the implications of North Carolina's HB92 0.1\% population equity constraint and how our algorithm compares, please refer to Appendix \ref{app:addressing-HB92}. This constraint is particularly challenging since other factors such as age distribution, voter registration, and turnout can affect the number of votes cast in a VTD, making population representation within the VTD incomplete. Further discussion on HB92's voter equity constraint is available in Appendix \ref{app:addressing-HB92}. While we could improve population equity by introducing more potential county splits, we believe our results demonstrate the effectiveness of our method. Our findings can serve as a starting point for other redistricting algorithms, or be fine-tuned post-hoc at the precinct level to satisfy strict population equity constraints, or serve as a baseline for comparing other redistricting algorithms.

The rest of the paper is outlined as follows.  Section \ref{sec:pre-processing} discusses preprocessing steps to prepare data for the algorithm, which is thoroughly outlined in Section \ref{sec:algorithm}.  Next, we generate a large sample of plausible redistricting plans according to the algorithm, and apply historical partisan votes for 2016 to the sampled plans as well as the historical North Carolina 2012, 2016 and 2020 plans, as well as the plan proposed by a bipartisan panel of retired NC judges, which has been deemed to be not gerrymandered (\cite{herschlag2020quantifying}). The results are discussed in Section \ref{sec:analysis}.  Lastly, an overall discussion is provided in Section \ref{sec:discussion}.  




\section{Pre-processing}\label{sec:pre-processing}

Looking at population density in North Carolina, we see that the county population varies from 4,407 to 919,628 among the 100 counties of North Carolina, ranging from $0.05\%$ to $9.64\%$ of the state population. In the pre-processing steps, we apply a splitting process to the three largest population counties.  This both helps to minimize county splits and allows for sufficient population share within each voting district (VTD). Splitting of counties can be controversial and can have negative ramifications to the fairness of a district plan, so we have tried to do the splitting in a reasonable manner.

\subsection{Splitting Counties with Large Population Size}
There are two main reasons that we choose to split these three counties. First and foremost, we want to minimize the number of county splits in order to prevent the potential of spreading voters of a particular type to multiple districts. Secondly, these 3 counties are the top 3 counties that have the largest population. Since our primary approach is to minimize the number of county splits, the top three counties are the best choice.

\begin{figure}[h!]
    \centering
    \includegraphics[scale=0.13]{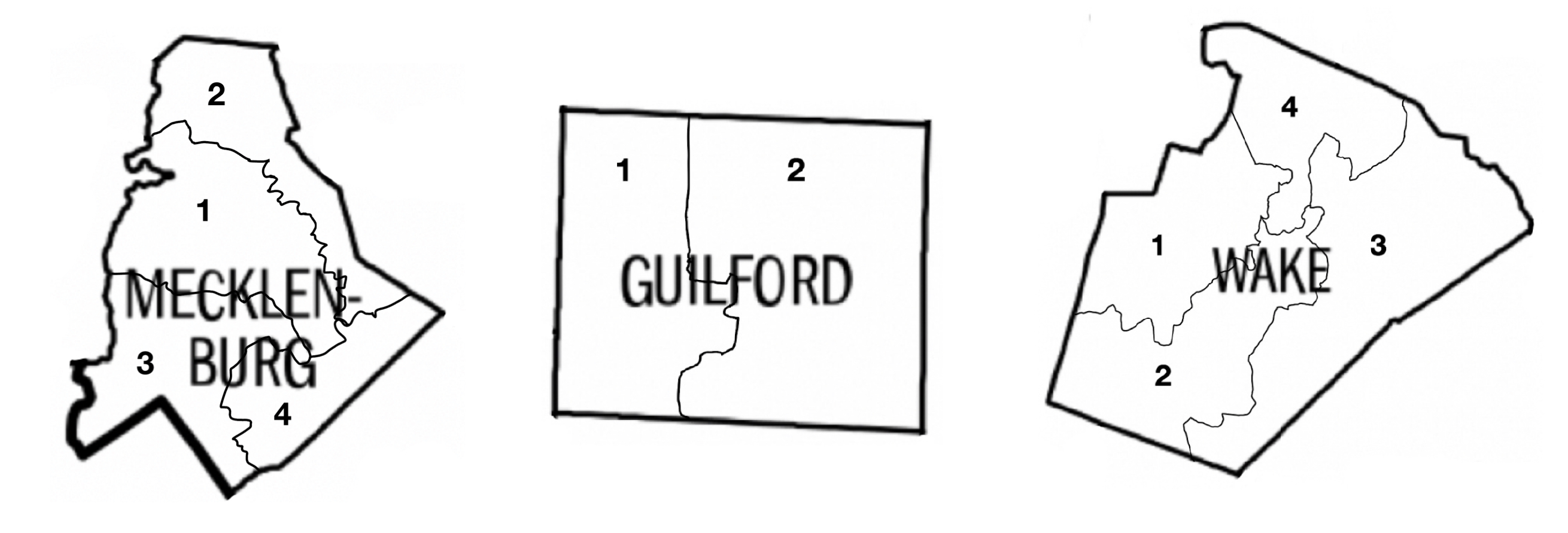}
    \caption{Our model's sub counties of Mecklenburg, Guilford, and Wake.}
    \label{fig:Our models' sub counties of Mecklenburg}
\end{figure}

Wake county and Mecklenburg county each have close to $10\%$ of the state population, and each will be split into four sub-counties. Guilford county has about $5\%$ of the state population, and it will be split into two sub-counties.  When we split the three largest counties, we make the cuts based on the population density map to ensure population equity. When possible, the split also maximizes potential shared borders with other \emph{generalized counties} (unsplit counties, or subcounties of split counties). Thus, we are increasing the number of basic territorial units from $100$ counties to $107$ generalized counties.  Figure \ref{fig:Our models' sub counties of Mecklenburg} shows how the splits were performed. 

For each county that gets split, we assume the number of votes for each sub-county is equal.
For example, Mecklenburg has a total of $9.64\%$ of the state population and $349,300$ votes, so each of the four sub-counties will have $2.41\%$ of the state population, which corresponds to $87,325$ votes in each sub-county. This is a reasonable assumption as the cut was made according to the population density map. Hence, boundaries of sub-counties can be ``dragged" in order to reach equal VTD population share, while maintaining the adjacency conditions. 

Then, we use North Carolina vote outcomes by precinct to approximate the party preference percentages in the sub-counties using the process described below.

\subsection{Votes Redistribution}

In the actual district plan used in the 2016 election, Mecklenburg, Guilford, and Wake counties each overlapped with two districts. For each of these counties, we use the actual votes distribution to guide our data to be used as input to our model, as described below.

\begin{figure}[h!]
    \centering
    \subfigure[\centering Mecklenburg 1]{{\includegraphics[scale=0.1877]{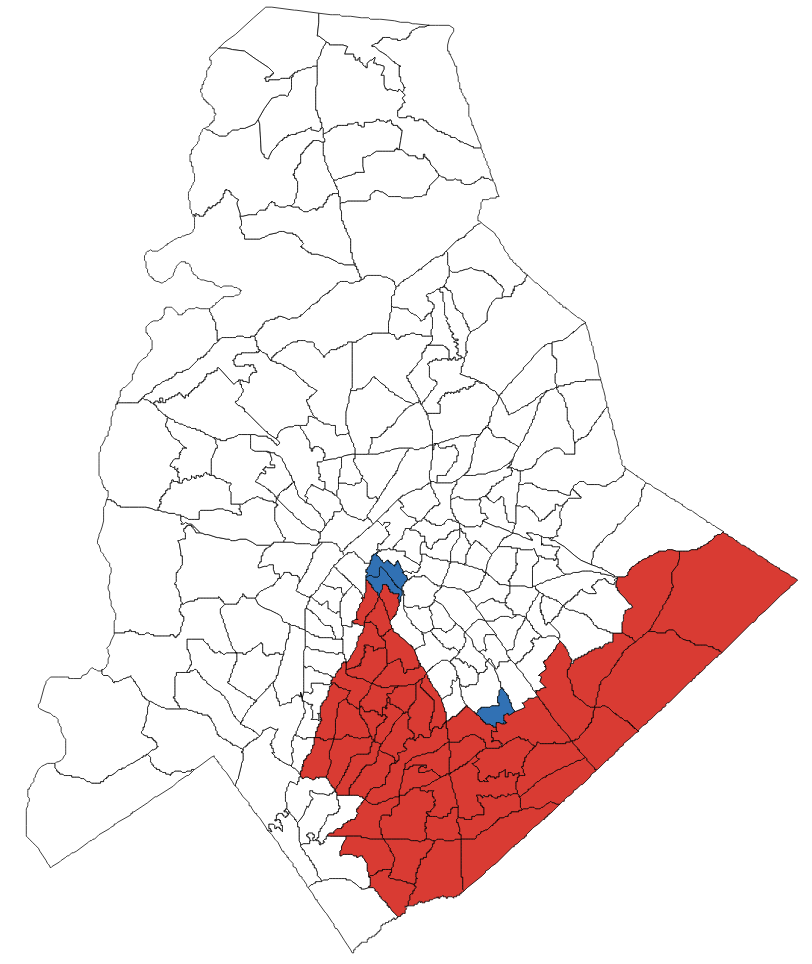} }}%
    \hspace{15mm}
    \subfigure[\centering Mecklenburg 2]{{\includegraphics[scale=0.1877]{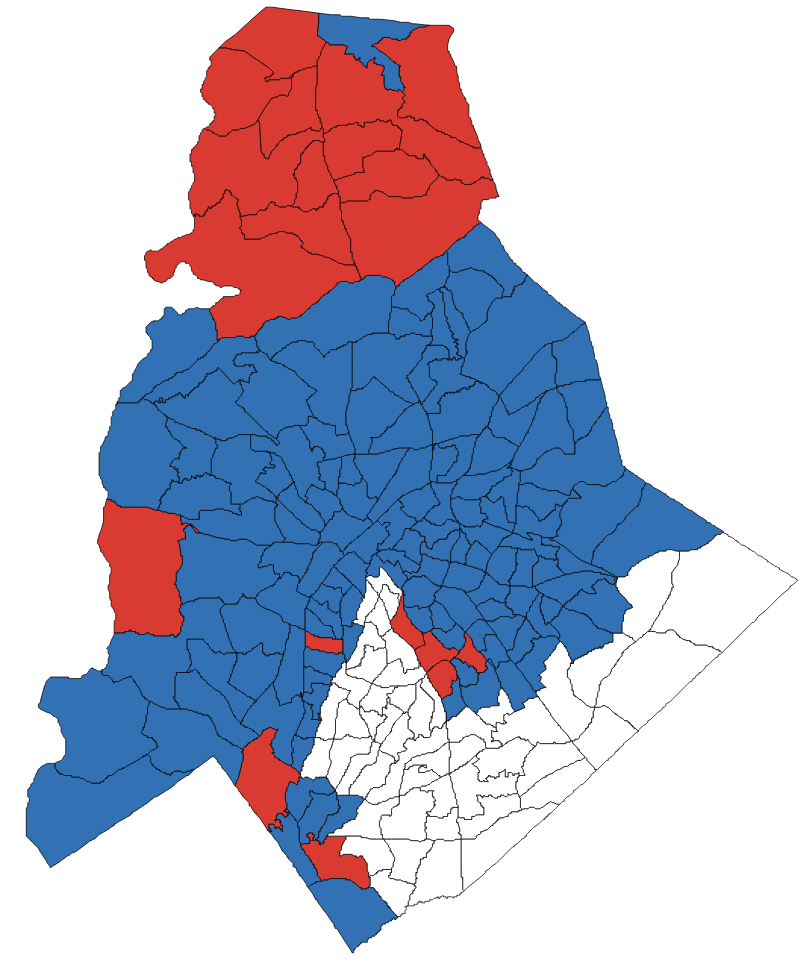} }}%
    \caption{2016 U.S. House of Representatives Election Results of Mecklenburg County by Precinct (Source: State of North Carolina Election Data, \url{https://er.ncsbe.gov/}).}
    \label{fig: 2016 U.S. House Representative Election Results of Mecklenburg County}
\end{figure}

\begin{table}[h!]
  \begin{tabular}{|c|r|r|r|r|}\hline
County	&	\multicolumn{2}{c|}{Democrat}			&	\multicolumn{2}{c|}{Republican}		\\ \hline
Mecklenburg 1	&	44,995	&	41.9\%	&	62,404	&	58.1\%	\\
(Belongs to District 9)	&		&		&		&		\\ \hline
Mecklenburg 2	&	234,115	&	67.0\%	&	115,185	&	33.0\%	\\
(Belongs to District 12)	&		&		&		&		\\ \hline \hline
\bf Total	&	279,110	&	61.1\%	&	177,589	&	38.9\%	\\ \hline
 \end{tabular}
\caption{Votes distribution of Mecklenburg County using actual 2016 district plan.}
\label{tab: Original Votes Distribution of Mecklenburg County}
\end{table} \vspace{13 mm}

\begin{figure}[h!]
    \centering
    \subfigure[\centering Guilford 1]{{\includegraphics[scale=0.1755]{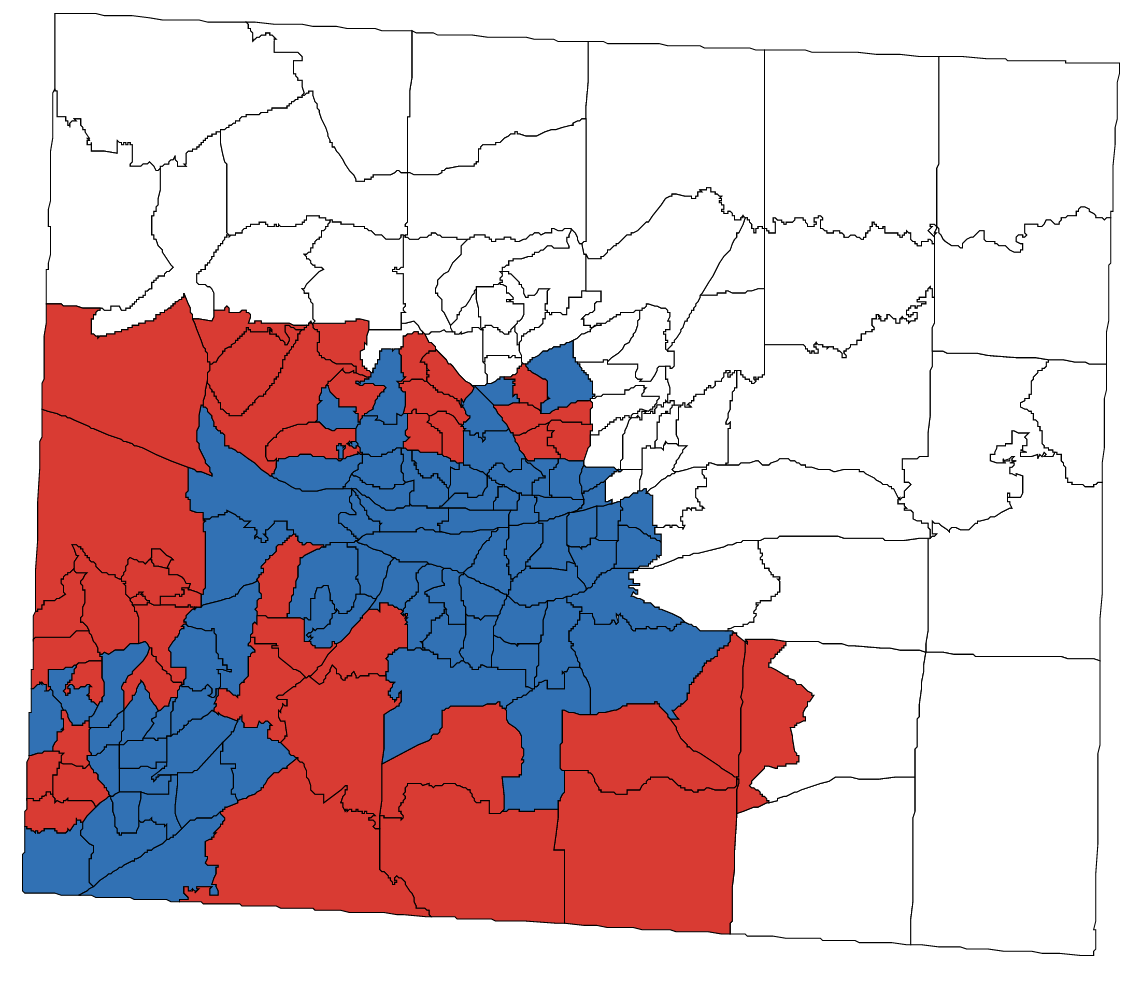} }}%
    \hspace{13 mm}
    \subfigure[\centering Guilford 2]{{\includegraphics[scale=0.1755]{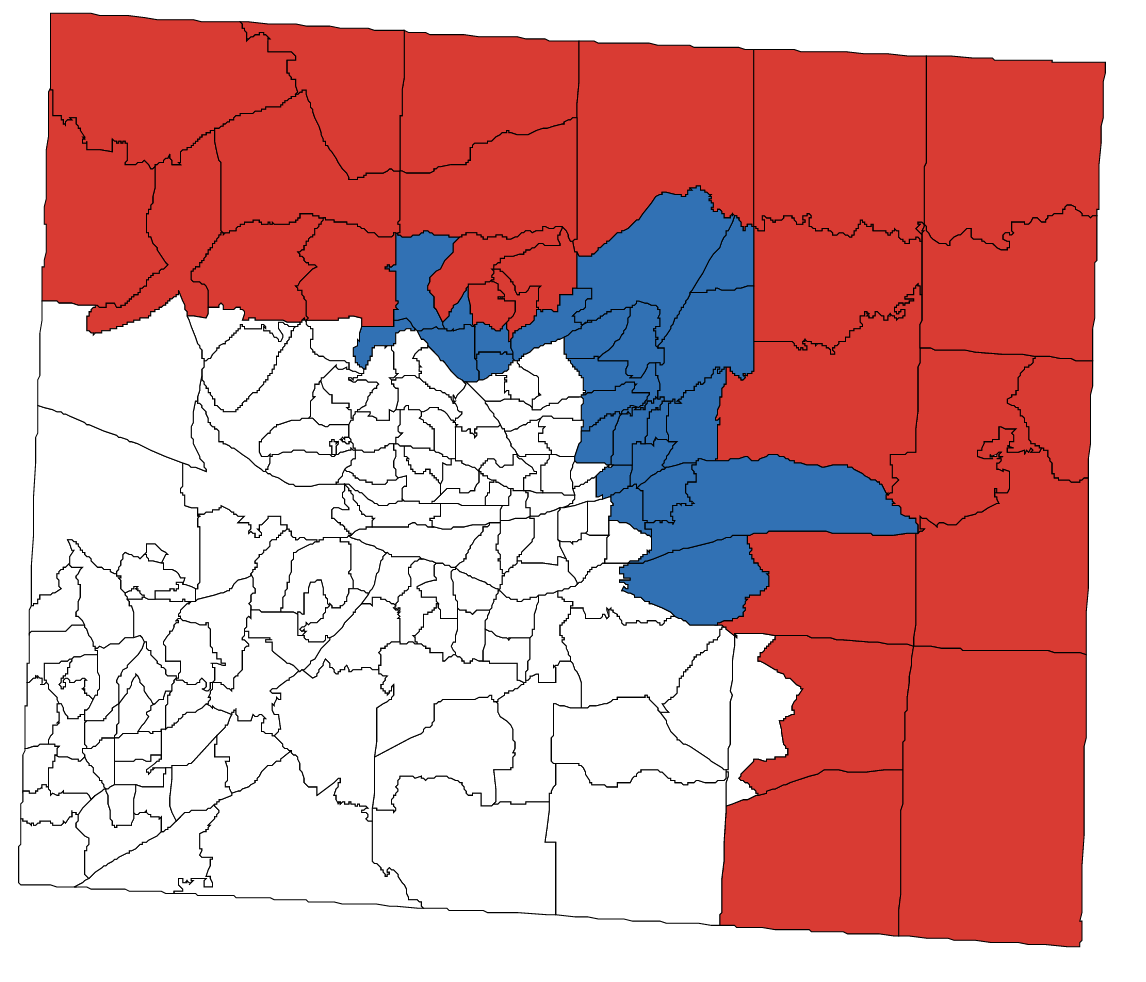} }}%
    \caption{2016 U.S. House of Representatives Election Results of Guilford County by Precinct (Source: State of North Carolina Election Data, \url{https://er.ncsbe.gov/}).}
    \label{fig: 2016 U.S. House Representative Election Results of Guilford County}
\end{figure}%

\begin{table}[h!]
  \begin{tabular}{|c|r|r|r|r|}\hline
County	&	\multicolumn{2}{c|}{Democrat}			&	\multicolumn{2}{c|}{Republican}			\\ \hline
Guilford 1	&	98,500	&	60.3\%	&	64,879	&	39.7\%	\\
(Belongs to District 13)	&		&		&		&		\\ \hline
Guilford 2	&	44,403	&	51.3\%	&	42,150	&	48.7\%	\\
(Belongs to District 6)	&		&		&		&		\\ \hline \hline
\bf Total	&	142,903	&	57.2\%	&	107,029	&	42.8\%	\\ \hline
 \end{tabular}
\caption{Votes Distribution of Guilford County using actual 2016 district plan.}
\label{tab: Original Votes Distribution of Guilford County}
\end{table}\vspace{13 mm}

\begin{figure}[!htb]
    \centering
    \subfigure[\centering Wake 1]{{\includegraphics[scale=0.1755]{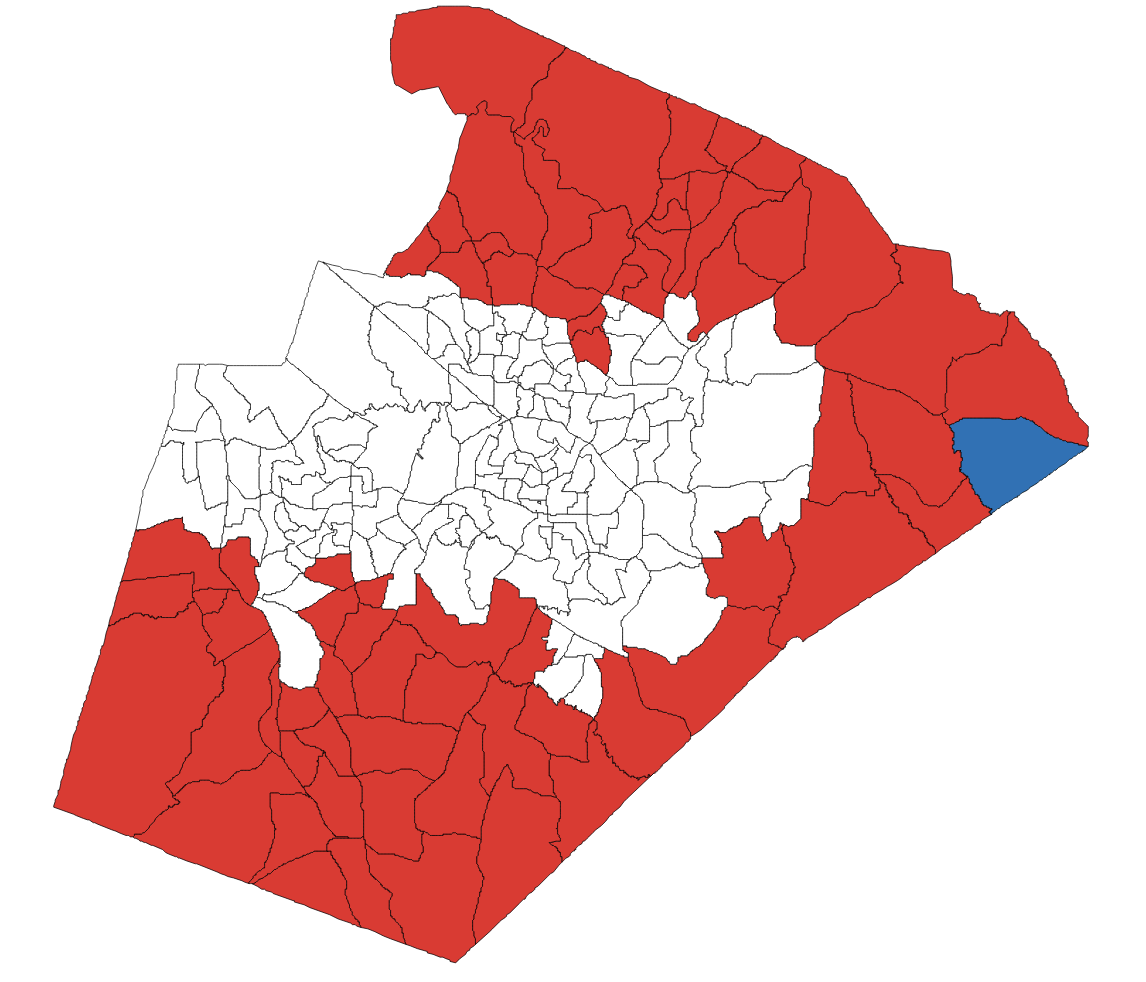} }}%
    \hspace{15mm}
    \subfigure[\centering Wake 2]{{\includegraphics[scale=0.1755]{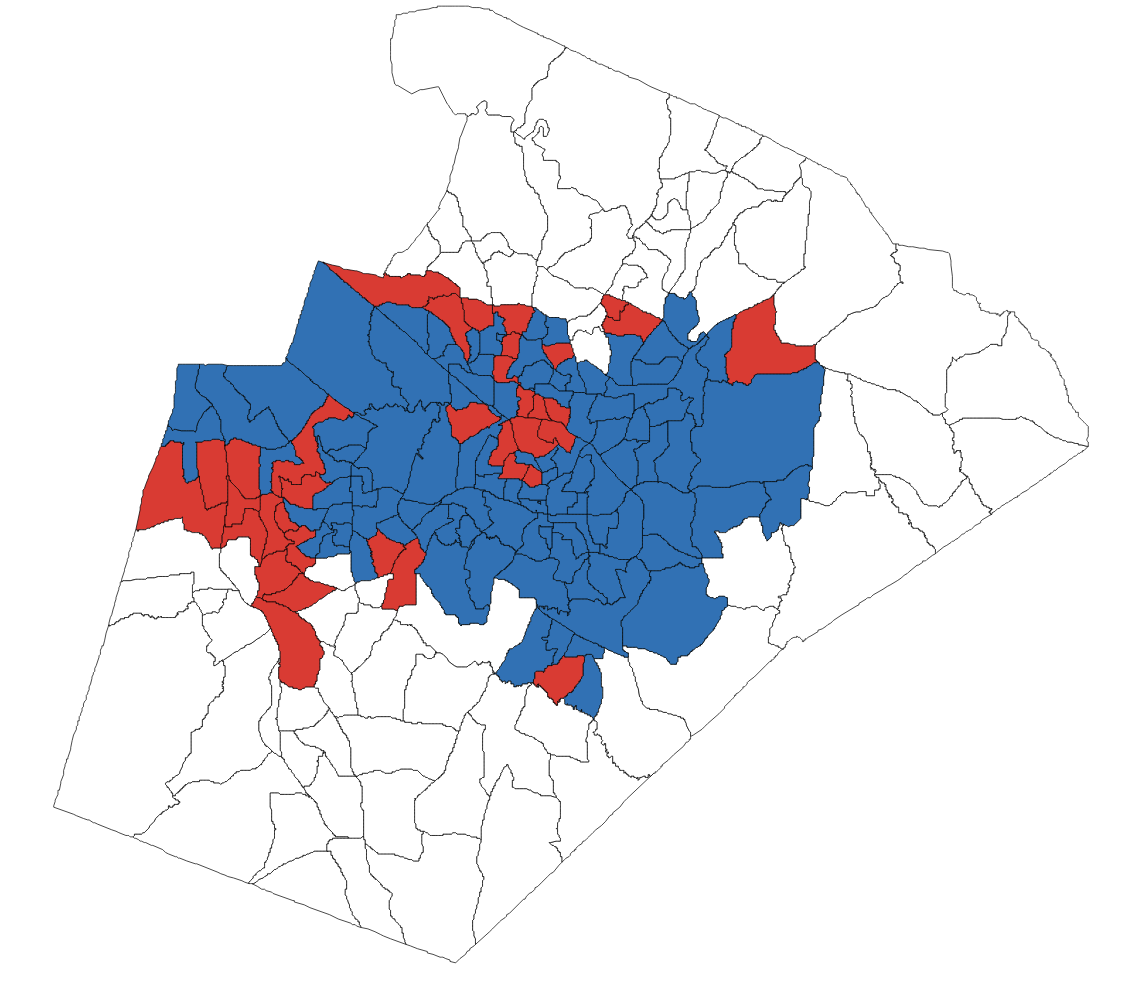} }}%
    \caption{2016 U.S. House of Representatives Election Results of Wake County by Precinct (Source: State of North Carolina Election Data, \url{https://er.ncsbe.gov/}).}
    \label{fig: 2016 U.S. House Representative Election Results of Wake County}
\end{figure}

\begin{table}[h!]
  \begin{tabular}{|c|r|r|r|r|}\hline
County	&	\multicolumn{2}{c|}{Democrat}			&	\multicolumn{2}{c|}{Republican}			\\ \hline
Wake 1	&	90,763	&	44.9\%	&	111,363	&	55.1\%	\\
(Belongs to District 2)	&		&		&		&		\\ \hline
Wake 2	&	206,172	&	66.2\%	&	105,327	&	33.8\%	\\
(Belongs to District 4)	&		&		&		&		\\ \hline \hline
Total	&	296,935	&	57.8\%	&	216,690	&	42.2\%	\\ \hline
 \end{tabular}
\caption{Votes Distribution of Wake County using actual 2016 district plan}
\label{tab: Original Votes Distribution of Wake County}
\end{table}

\begin{table}[ht]
   \begin{tabular}{|c|r|r|r|r|}\hline
County	&	\multicolumn{2}{c|}{Democrat}			&	\multicolumn{2}{c|}{Republican}				\\ \hline
Mecklenburg 1	&	88,257	&$	77.3	\%$&	25,917	&$	22.7	\%$	\\ \hline
Mecklenburg 2	&	73,529	&$	64.4	\%$&	40,646	&$	35.6	\%$	\\ \hline
Mecklenburg 3	&	66,745	&$	58.5	\%$&	47,430	&$	41.5	\%$	\\ \hline
Mecklenburg 4	&	50,580	&$	44.3	\%$&	63,595	&$	55.7	\%$	\\ \hline
   \end{tabular}
\caption{Votes Distribution of Mecklenburg County using our model's sub counties.} 
\label{tab: Votes Redistribution of Mecklenburg County}
\end{table}
\begin{table}[ht]
   \begin{tabular}{|c|r|r|r|r|}\hline
County	&	\multicolumn{2}{c|}{Democrat}			&	\multicolumn{2}{c|}{Republican}				\\ \hline
Guilford 1	&	69,143	&$	55.3	\%$&	55,823	&$	44.7	\%$	\\ \hline
Guilford 2	&	73,760	&$	59.0	\%$&	51,206	&$	41.0	\%$	\\ \hline
   \end{tabular}
\caption{Votes Distribution of Guilford County using our model's sub-counties.} 
\label{tab: Votes Redistribution of Guilford County}
\end{table}
\begin{table}[ht]
   \begin{tabular}{|c|r|r|r|r|}\hline
County	&	\multicolumn{2}{c|}{Democrat}	&	\multicolumn{2}{c|}{Republican}	\\ \hline
Wake 1	&	75,323	&$	58.9	\%$&	52,621	&$	41.1	\%$	\\ \hline
Wake 2	&	77,867	&$	61.6	\%$&	48,538	&$	38.4	\%$	\\ \hline
Wake 3	&	78,125	&$	59.3	\%$&	53,606	&$	40.7	\%$	\\ \hline
Wake 4	&	65,620	&$	51.4	\%$&	61,925	&$	48.6	\%$	\\ \hline
   \end{tabular}
\caption{Votes Distribution of Wake County using our model's sub counties.} 
\label{tab: Votes Redistribution of Wake County}
\end{table}

\subsubsection{Mecklenburg County}

In our splitting of Mecklenburg into four sub-counties, Mecklenburg 4 is similar to the Mecklenburg county area that belongs to the actual District~9. Mecklenburg 1, 2, and 3 combined are similar to the Mecklenburg county area that belongs to the actual District~12. Figure~\ref{fig: 2016 U.S. House Representative Election Results of Mecklenburg County} shows visually the precinct data for the portions of Mecklenburg belonging to the actual District~9 and District~12, respectively. We approximate the votes for each sub-county accordingly. 
The actual vote distribution for Mecklenburg County is shown in Table~\ref{tab: Original Votes Distribution of Mecklenburg County}.
Our splitting of Mecklenburg is shown in Figure~\ref{fig:Our models' sub counties of Mecklenburg}.
The final vote redistribution for our model's sub-counties of Mecklenburg is shown in Table~\ref{tab: Votes Redistribution of Mecklenburg County}.  
\subsubsection{Guilford County}

Similar to Mecklenburg county, the votes of Guilford county were also counted for the sub-counties of Guilford county that we created. See Figure~\ref{fig: 2016 U.S. House Representative Election Results of Guilford County} and Table~\ref{tab: Original Votes Distribution of Guilford County} for the actual precinct data and vote distribution. See Figure~\ref{fig:Our models' sub counties of Mecklenburg} and Table~\ref{tab: Votes Redistribution of Guilford County} for our splitting of Guilford and the resulting vote distribution.

\subsubsection{Wake County}

Similar to Mecklenburg county, the votes of Wake county were also counted for the sub-counties of Wake county that we created. See Figure~\ref{fig: 2016 U.S. House Representative Election Results of Wake County} and Table~\ref{tab: Original Votes Distribution of Wake County} for the actual precinct data and vote distribution. See Figure~\ref{fig:Our models' sub counties of Mecklenburg} and Table~\ref{tab: Votes Redistribution of Wake County} for our splitting of Wake and the resulting vote distribution.

\subsection{Additional Pre-processing}

In our modeling, adjacency data for North Carolina counties from the U.S. Census is used. However, there are some ambiguous cases due to weak borders. Take Halifax county and Franklin county as an example; they do share a ``border" since there are certain parts of their land that connect, but these connections are extremely weak, meaning they barely connect to each other.

Since we will cluster multiple counties into one district, clustering two counties that share an extremely weak connection can result in a district that is not sufficiently compact. Thus, we decide that such pairs of counties are not viewed as adjacent. The full list of our model's non-adjacent county pairs that are different from the U.S. Census data is shown below.

However, this does not mean that these pairs of non-adjacent counties do not have a chance to be included in the same district. Instead, if their common adjacent county gets selected in a district, the pair will have a chance to be included in the later selection process.

\begin{paracol}{2}
\begin{itemize}
\item	Halifax	-	Franklin
\item	Henderson	-	Haywood
\item	Johnston	-	Franklin
\item	Lincoln	-	Burke
\item	Montgomery	-	Anson
\item	Moore	-	Cumberland
\item	Richmond	-	Hoke
\item	Rockingham	-	Alamance
\item	Rockingham	-	Forsyth

\end{itemize}

\switchcolumn

\begin{itemize}
\item	Rowan	-	Montgomery
\item	Scotland	-	Moore
\item	Stanly	-	Davidson
\item	Stanly	-	Richmond
\item	Stokes	-	Guilford
\item	Surry	-	Forsyth
\item	Wake	-	Nash
\item	Warren	-	Nash
\end{itemize}
\end{paracol}

\begin{figure}[h!]
    \centering
    \includegraphics[scale=0.17]{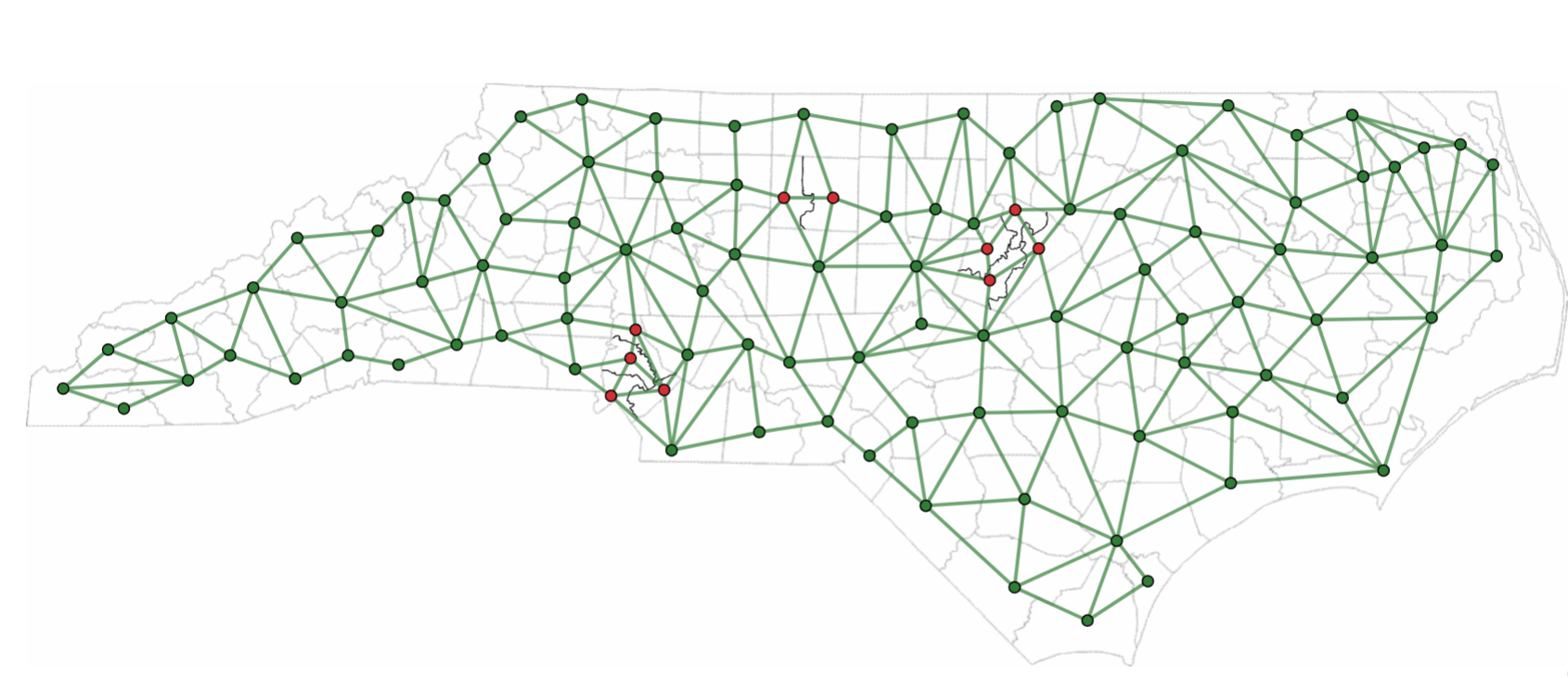}
    \caption{Adjacency network of 107 generalized counties for North Carolina.}
    \label{fig:NC_network by MATLAB}
\end{figure}

As a result of the pre-processing work, we obtain an adjacency network with $107$ nodes as shown in Figure~\ref{fig:NC_network by MATLAB}, representing the $107$ generalized counties. An edge between two nodes indicates that the corresponding generalized counties are adjacent. The red nodes are the sub-counties of the three split counties. The associated adjacency matrix is used extensively in our clustering algorithm for forming district plans.

In our algorithm, when it selects the third adjacent county, the county that shares more borders with the selected counties (which can be interpreted as the one that has higher \% border shared) will have higher chances than other counties to be selected.
\section{Algorithm for Redistricting of North Carolina}\label{sec:algorithm}
\renewcommand{\labelenumii}{\arabic{enumi}.\arabic{enumii}}
\renewcommand{\labelenumiii}{\arabic{enumi}.\arabic{enumii}.\arabic{enumiii}}
\renewcommand{\labelenumiv}{\arabic{enumi}.\arabic{enumii}.\arabic{enumiii}.\arabic{enumiv}}

We now describe our clustering algorithm for creating district plans for North Carolina, composed of thirteen districts clustered from $107$ generalized counties. Each district to be formed is given a distinct index from $1$, $2$, \ldots $13$. Each generalized county is given a distinct index from $1$, $2$, \ldots $107$.

Our algorithm pre-populates the thirteen districts using thirteen ``seeds.'' The seeds are manually chosen with intention. More details on pre-selected counties are discussed in section \ref{sec: Selecting Criteria of Pre-selected Counnties}. By the design of this algorithm, the pre-selected seed for each district tends to have a large population size among its neighboring counties. This does potentially rule out the possibility that a district only contains generalized counties with relatively small population sizes. Analysis of results obtained from various sets of seeds indicates strong robustness to our results; many sets of seeds do not yield viable results, and those sets of seeds that do yield viable results give similar conclusions to those we will present in our analysis.

After the initial ``seeding'', our algorithm then begins adding to each district large populations of generalized counties, and then later on adding smaller generalized counties. 
We use an analogy for our algorithm of adding to a bottle rocks, then pebbles, then sand. 

\subsection{Clustering Algorithm}

Here is a description of our algorithm:

\begin{enumerate}
    \item Create the $107 \times 107$ generalized county adjacency matrix $A$.
    \item Create a $60 \times 13$ matrix, $Sol$, in which each column corresponds to a distinct district.
    \item ($i==1$) Input 1st generalized county for each of the 13 districts:
            \begin{enumerate}
            \item $Sol(1,j)$ is filled with pre-selected generalized county indices.
        \end{enumerate}
    \item ($i==2$) Choose 2nd generalized county for each of the 13 districts: 
        \begin{enumerate}
            \item For each column $j$, find all the un-selected indices of adjacent generalized counties of the 1st generalized county  in column $j$.
            \item Among these indices, pick the generalized county index that has the largest population to become the second generalized county of the district, $Sol(2,j)$. In the case of a tie due to county splits, the tie is broken randomly.
        \end{enumerate}
     \item ($i==3$) Choose 3rd generalized county for each of the 13 districts: 
        \begin{enumerate}
            \item For each column $j$, find all the un-selected indices of adjacent generalized counties of the 1st generalized county and all the un-selected indices of adjacent generalized counties of the 2nd generalized county in column $j$.
            \item Among these indices, pick the generalized county index that has the largest population to become the 3rd generalized county of the district, $Sol(3,j)$. In the case of a tie due to county splits, the tie is broken randomly.
        \end{enumerate}
     \item ($i==4,5$) Choose additional generalized counties for the $13$ districts: 
        \begin{enumerate}
            \item  Case 1: For each column $j$, if the developing district's population is $>6\%$ of the state population, then $Sol(i,j)=0$. 
            \item  Case 2: For each column $j$, if the developing district's population is $\leq6\%$ of the state population, then find all the un-selected indices of adjacent generalized counties for the existing generalized counties in column $j$ and randomly pick one of them to be included in $Sol(i,j)$.
        \end{enumerate}
    \item ($6 \leq i$) Choose additional generalized counties for each of the 13 districts: 
        \begin{enumerate}
        \item  Case 1: For each column $j$, if the district's population is $>\Bigg(6+ \Big\lfloor \dfrac{i-3}{3}\Big\rfloor \cdot 0.45\Bigg)\%$ of the state population, then $Sol(i,j)=0$.
        \item  Case 2: For each column $j$, if the district's population is $\leq\Bigg(6+ \Big\lfloor \dfrac{i-3}{3}\Big\rfloor \cdot 0.45\Bigg)\%$ of the state population, then find all the adjacent generalized counties for the existing generalized counties in column $j$ and randomly pick one of them to be included in $Sol(i,j)$.
        \item If all generalized county indices have been selected, then the matrix $Sol$ is complete. Break from the current ``for" loop. Otherwise, repeat step 7 until all generalized county indices have been selected.
        \end{enumerate}
\end{enumerate}

\subsection{Algorithm Notes}

\begin{figure}[h!]
\centering
        \includegraphics[scale=0.1]{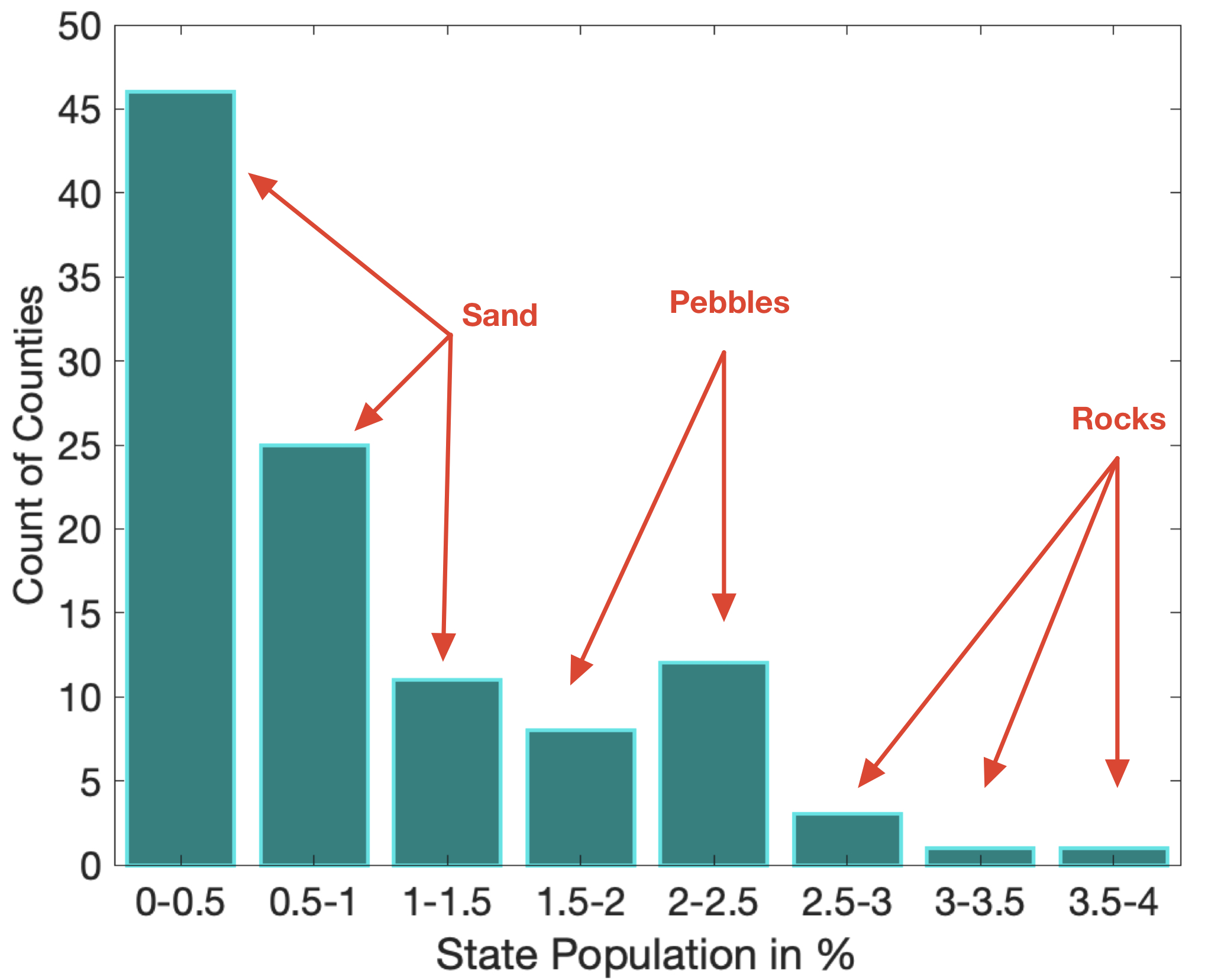}
        \caption{Summary of Population Size Among Generalized Counties in North Carolina (A Total of 107 Counties After Splitting)}
        \label{fig: Rock-pebbles-sand}
\end{figure}

\begin{itemize}
    
    \item We pre-allocated $60$ rows to $Sol$ to provide a sufficient number of iterations for simulations to settle in to completion of the matrix $Sol$.
    
    \item In the randomized process (starting when $i\geq 4$) of selecting generalized counties from the adjacency list, generalized counties that are adjacent to multiple generalized counties in a developing district will be included in the list of adjacent, unselected indices  multiple times, resulting in a higher probability that they will get selected. This gives preference to the development of districts that are relatively compact.
    
    \item The first three generalized counties of each district have large population size. The purpose of this is to avoid the possibility of grouping many generalized counties with small populations into one district in the beginning, which is likely to lead to a larger population inequity among the districts. The only randomization in the second and third generalized counties of each district is due to the breaking of population ties caused by county splits.
    
    \item Substantial randomization on selecting generalized counties starts at the 4th row for each district. Population percentage caps are used to let districts that have relatively small populations ``catch up'' before letting other districts that have relatively large populations continue including more generalized counties into their districts.
    
    \item Starting with the 6th row, we use a different population percentage cap, $\Bigg(6+ \Big\lfloor \dfrac{i-3}{3}\Big\rfloor \cdot 0.45\Bigg)\%$, and the cap jumps at a constant rate every 3 rows thereafter. The percentage $0.45\%$ is derived from Figure~\ref{fig: Rock-pebbles-sand}. There is a total of $107$ generalized counties, and close to half of these have population under $0.5\%$ of the state population. While testing different combinations of numbers for a population cap, we found that usage of $0.45\%$ results in a relatively low standard deviation in district population in resulting district plans. Usage of the floor function gives more room for generalized counties that have similar population sizes to be utilized before the cap is increased again.
    
    \begin{itemize}
    
        \item If the percentage that was chosen is too small, say $0.1\%$ instead of $0.45\%$, then many more iterations are required for the process to terminate, and the resulting district plans tend to have
        large standard deviations in district population.
    
        \item If the number that was chosen is too large, say $1\%$ instead of $0.45\%$, although the algorithm terminates rapidly, the resulting district plans again tend to have
        large standard deviations in district population.
    
    \end{itemize}

\end{itemize}

\subsection{Selecting Criteria of Pre-selected Counties}\label{sec: Selecting Criteria of Pre-selected Counnties}
In Figure \ref{fig: Rock-pebbles-sand}, we categorized the 107 generalized counties by population size into 3 categories: rocks, pebbles, and sand. By the design of this algorithm, the pre-selected first counties of each district are either rocks or pebbles, or they have a very large population or are densely populated compared to their neighboring counties. Readers can verify that from Table \ref{tab: Population of 107 Generalized Counties}. The counties in red in Table \ref{tab: Population of 107 Generalized Counties} are the pre-selected counties.

\begin{figure}[htbp!]
    \centering
    \includegraphics[scale=0.25]{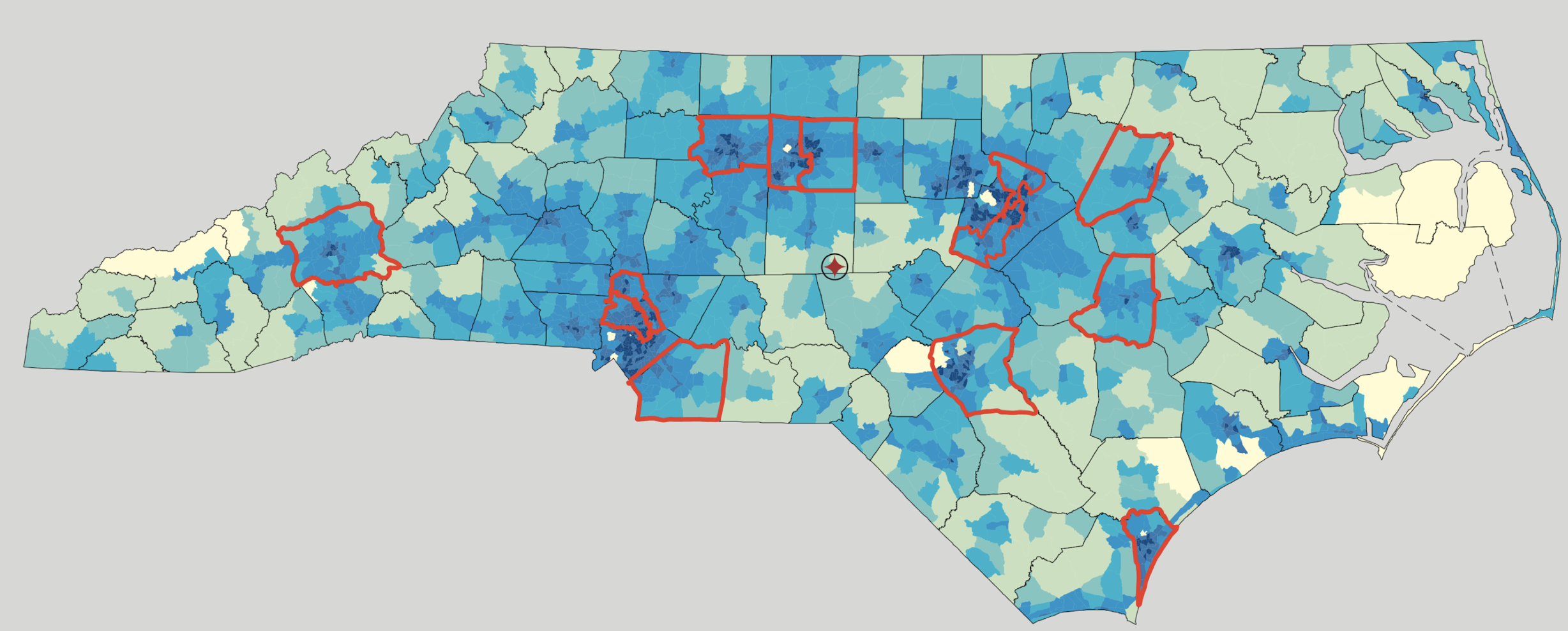}
    \caption{Pre-selected Counties on North Carolina Population Density map of 2010 Census (Source: United States Census Bureau, \url{https://www2.census.gov/geo/maps/dc10_thematic/2010_Profile/2010_Profile_Map_North_Carolina.pdf})}
    \label{fig:Pre-selected Counites on Population Density map}
\end{figure}

Although Wayne and Nash counties do not belong to Pebbles, they have the largest population or population density around their neighboring counties as shown in Figure \ref{fig:Pre-selected Counites on Population Density map}.
We are aware that there are other combinations of such criteria that exist, and the reasons for selecting the particular set of the pre-selected counties in this study are based on the combination of the following:

\begin{itemize}
    \item Location of the counties 
    \begin{itemize}
        \item high-density areas: counties have more adjacent counties 
        \item not high-density areas: has the largest population and/or population density around its neighboring counties. 
    \end{itemize}
    \item Population density around its neighboring counties and population size.
    \item We have tested numerous sets of pre-selected counties to see which ones can provide more plausible district maps while maintaining strong population equity.
\end{itemize}

With careful consideration, we chose to show this best version of the result from these pre-selected counties as shown in Figure \ref{fig:Pre-selected Counites on Population Density map} and Table \ref{tab: Population of 107 Generalized Counties}.

If the state would like certain counties to be in separate districts, they can use our model and algorithm with modified pre-selected counties.

\begin{table}[h!]
    \centering
    \resizebox{\columnwidth}{!}{\begin{tabular}{|l|l|c|c|l|l|c|}\hline
Rank	&	County	&	\% of State Population		&&	Rank	&	County	&	\% of State Population	\%\\\hline
\textcolor{red}{1}	&	\textcolor{red}{Forsyth}	&	\textcolor{red}{3.68	\%}	&&	55	&	Columbus	&	0.61	\%\\ 
\textcolor{red}{2}	&	\textcolor{red}{Cumberland}	&	\textcolor{red}{3.35	\%}	&&	56	&	Lee	&	0.61	\%\\ 
3	&	Durham	&	2.81	\%	&&	57	&	Edgecombe	&	0.59	\%\\
\textcolor{red}{4}	&	\textcolor{red}{Guilford 1}	&	\textcolor{red}{2.56	\%}	&&	58	&	Halifax	&	0.57	\%\\
\textcolor{red}{5}	&	\textcolor{red}{Guilford 2}	&	\textcolor{red}{2.56	\%}	&&	59	&	Pender	&	0.55	\%\\
\textcolor{red}{6}	&	\textcolor{red}{Buncombe}	&	\textcolor{red}{2.50	\%}&&	60	&	Watauga	&	0.54	\%\\
\textcolor{red}{7}	&	\textcolor{red}{Mecklenburg 1}	&	\textcolor{red}{2.41	\%}	&&	61	&	Beaufort	&	0.50	\%\\
\textcolor{red}{8}	&	\textcolor{red}{Mecklenburg 2}	&	\textcolor{red}{2.41	\%}	&&	62	&	Stokes	&	0.50	\%\\
9	&	Mecklenburg 3	&	2.41	\%	&&	63	&	Hoke	&	0.49	\%\\
10	&	Mecklenburg 4	&	2.41	\%	&&	64	&	Richmond	&	0.49	\%\\
\textcolor{red}{11}	&	\textcolor{red}{Wake 4}	&	\textcolor{red}{2.36	\%}	&&	65	&	Vance	&	0.48	\%\\
12	&	Wake 1	&	2.36	\%	&&	66	&	McDowell	&	0.47	\%\\
\textcolor{red}{13}	&	\textcolor{red}{Wake 2}	&	\textcolor{red}{2.36	\%}&&	67	&	Davie	&	0.43	\%\\
14	&	Wake 3	&	2.36	\%	&&	68	&	Pasquotank	&	0.43	\%\\
15	&	Gaston	&	2.16	\%	&&	69	&	Jackson	&	0.42	\%\\
\textcolor{red}{16}	&	\textcolor{red}{New Hanover}	&	\textcolor{red}{2.13	\%}	&&	70	&	Person	&	0.41	\%\\
\textcolor{red}{17}	&	\textcolor{red}{Union}	&	\textcolor{red}{2.11	\%}	&&	71	&	Yadkin	&	0.40	\%\\
18	&	Cabarrus	&	1.87	\%	&&	72	&	Alexander	&	0.39	\%\\
19	&	Onslow	&	1.86	\%	&&	73	&	Scotland	&	0.38	\%\\
20	&	Johnston	&	1.77	\%	&&	74	&	Bladen	&	0.37	\%\\
21	&	Pitt	&	1.76	\%	&&	75	&	Macon	&	0.36	\%\\
22	&	Davidson	&	1.71	\%	&&	76	&	Dare	&	0.36	\%\\
23	&	Iredell	&	1.67	\%	&&	77	&	Transylvania	&	0.35	\%\\
24	&	Catawba	&	1.62	\%	&&	78	&	Montgomery	&	0.29	\%\\
25	&	Alamance	&	1.58	\%	&&	79	&	Cherokee	&	0.29	\%\\
26	&	Randolph	&	1.49	\%	&&	80	&	Ashe	&	0.29	\%\\
27	&	Rowan	&	1.45	\%	&&	81	&	Anson	&	0.28	\%\\
28	&	Robeson	&	1.41	\%	&&	82	&	Hertford	&	0.26	\%\\
29	&	Orange	&	1.40	\%	&&	83	&	Martin	&	0.26	\%\\
\textcolor{red}{30}	&	\textcolor{red}{Wayne}	&	\textcolor{red}{1.29	\%}	&&	84	&	Caswell	&	0.25	\%\\
31	&	Harnett	&	1.20	\%	&&	85	&	Currituck	&	0.25	\%\\
32	&	Brunswick	&	1.13	\%	&&	86	&	Northampton	&	0.23	\%\\
33	&	Henderson	&	1.12	\%	&&	87	&	Greene	&	0.22	\%\\
34	&	Craven	&	1.09	\%	&&	88	&	Bertie	&	0.22	\%\\
35	&	Cleveland	&	1.03	\%	&&	89	&	Warren	&	0.22	\%\\
\textcolor{red}{36}	&	\textcolor{red}{Nash}	&	\textcolor{red}{1.01	\%}	&&	90	&	Madison	&	0.22	\%\\
37	&	Rockingham	&	0.98	\%	&&	91	&	Polk	&	0.22	\%\\
38	&	Burke	&	0.95	\%	&&	92	&	Yancey	&	0.19	\%\\
39	&	Moore	&	0.93	\%	&&	93	&	Avery	&	0.19	\%\\
40	&	Caldwell	&	0.87	\%	&&	94	&	Mitchell	&	0.16	\%\\
41	&	Wilson	&	0.85	\%	&&	95	&	Chowan	&	0.16	\%\\
42	&	Lincoln	&	0.82	\%	&&	96	&	Swain	&	0.15	\%\\
43	&	Surry	&	0.77	\%	&&	97	&	Perquimans	&	0.14	\%\\
44	&	Wilkes	&	0.73	\%	&&	98	&	Washington	&	0.14	\%\\
45	&	Rutherford	&	0.71	\%	&&	99	&	Pamlico	&	0.14	\%\\
46	&	Carteret	&	0.70	\%	&&	100	&	Gates	&	0.13	\%\\
47	&	Chatham	&	0.67	\%	&&	101	&	Alleghany	&	0.12	\%\\
48	&	Sampson	&	0.67	\%	&&	102	&	Clay	&	0.11	\%\\
49	&	Franklin	&	0.64	\%	&&	103	&	Jones	&	0.11	\%\\
50	&	Stanly	&	0.64	\%	&&	104	&	Camden	&	0.10	\%\\
51	&	Granville	&	0.63	\%	&&	105	&	Graham	&	0.09	\%\\
52	&	Lenoir	&	0.62	\%	&&	106	&	Hyde	&	0.06	\%\\
53	&	Haywood	&	0.62	\%	&&	107	&	Tyrrell	&	0.05	\%\\
54	&	Duplin	&	0.61	\%	&&		&		&		\\\hline
    \end{tabular}}
        \caption{Population of 107 Generalized Counties (Ranked by Population Size)}
    \label{tab: Population of 107 Generalized Counties}
\end{table}


\section{Ensemble of District Plans}\label{sec:ensemble}

\begin{figure}[h!]
    \centering
    \includegraphics[scale=0.55]{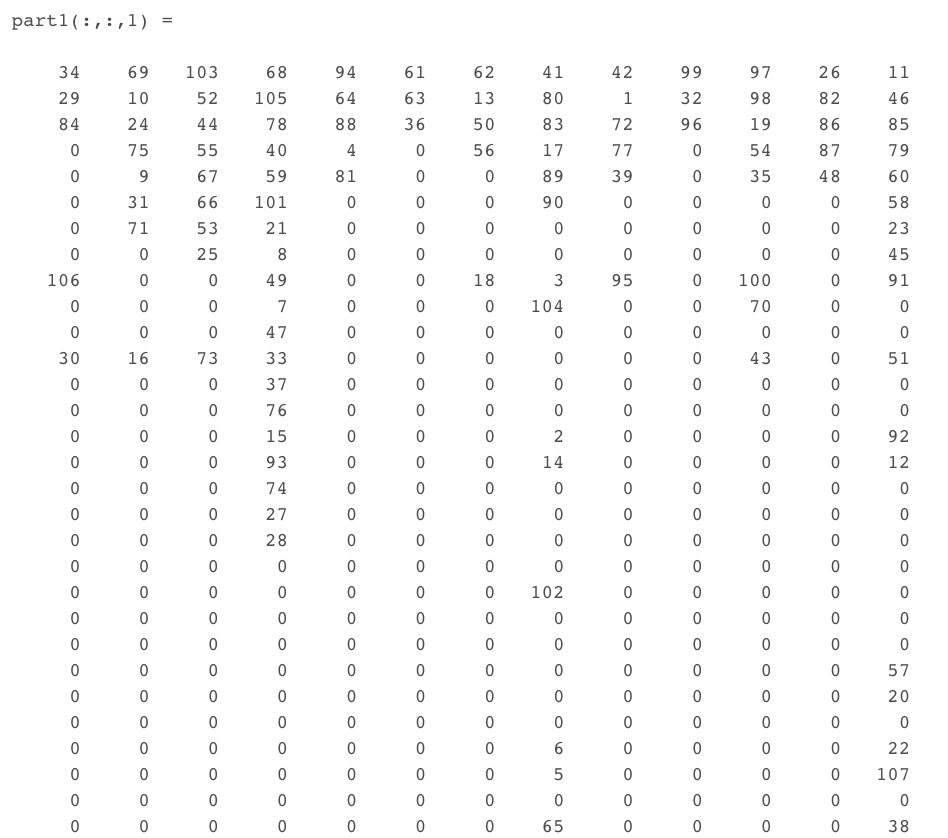}
    \caption{Sample Result 1 of a potential redistricting plan for North Carolina from MATLAB. Each entry (except 0) represents a generalized county and each column represents a district. For display purposes, we have truncated the output not to show rows consisting of all zeros.}
    \label{fig:Sample Result 1 on Potential Redistricting plan for North Carolina from MATLAB}
    \includegraphics[scale=0.55]{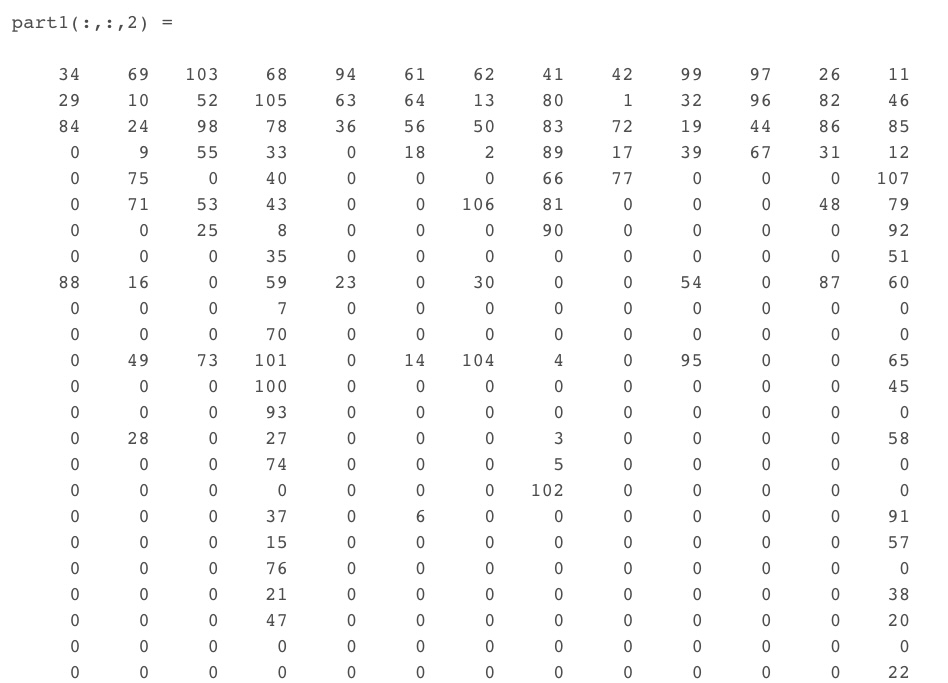}
    \caption{Sample Result 2 of a potential redistricting plan for North Carolina from MATLAB.}
    \label{fig:Sample Result 2 on Potential Redistricting plan for North Carolina from MATLAB}
\end{figure}

\begin{table}[htbp!]
    \resizebox{\textwidth}{!}{
    \centering
\begin{tabular}{|c|c|c|c|c|c|c|c|c|c|c|c|c|c|}\hline								
District		&	1	&	2	&	3	&	4	&	5	&	6	&	7	&	8	&	9	&	10	&	11	&	12	&	13	\\ \hline
	\multirow{18}{*}{County Indices}	&   34	&	69	&	103	&	68	&	94	&	61	&	62	&	41	&	42	&	99	&	97	&	26	& 	11	\\
		&	29	&	10	&	52	&	105	&	63	&	64	&	13	&	80	&	1	&	32	&	98	&	82	&	46	\\
		&	84	&	24	&	44	&	78	&	36	&	56	&	50	&	83	&	72	&	96	&	19	&	86	&	85	\\
		&	106	&	75	&	55	&	7	&	23	&	18	&	2	&	17	&	77	&		&	67	&	87	&	58	\\
		&	30	&	9	&	54	&	73	&	4	&	12	&	104	&	89	&	39	&		&	35	&	48	&	45	\\
		&		&	31	&	40	&	33	&		&		&	14	&	90	&	95	&		&		&	81	&	91	\\
		&		&	71	&	25	&	49	&		&		&		&	3	&	100	&		&		&	66	&	79	\\
		&		&	53	&	16	&	101	&		&		&		&	5	&	70	&		&		&	88	&	92	\\
		&		&		&		&	28	&		&		&		&	102	&	47	&		&		&		&	107	\\
		&		&		&		&	59	&		&		&		&	6	&	37	&		&		&		&	51	\\
		&		&		&		&	43	&		&		&		&	65	&		&		&		&		&	60	\\
		&		&		&		&	8	&		&		&		&		&		&		&		&		&	38	\\
		&		&		&		&	93	&		&		&		&		&		&		&		&		&	20	\\
		&		&		&		&	21	&		&		&		&		&		&		&		&		&	57	\\
		&		&		&		&	27	&		&		&		&		&		&		&		&		&	22	\\
		&		&		&		&	76	&		&		&		&		&		&		&		&		&		\\
		&		&		&		&	74	&		&		&		&		&		&		&		&		&		\\
		&		&		&		&	15	&		&		&		&		&		&		&		&		&		\\ \hline
State Population (\%)		&$	7.7	\%$&$	7.4	\%$&$	7.5	\%$&$	7.6	\%$&$	8.0	\%$&$	8.2	\%$&$	7.9	\%$&$	7.8	\%$&$	7.9	\%$&$	7.5	\%$&$	7.0	\%$&$	7.7	\%$&$	7.8	\%$ \\ \hline
    \end{tabular}}
    \caption{Best solution (in terms of population standard deviation) of the simulation for redistricting of North Carolina with population distribution.}
    \label{tab: Best Solution of the Simulation for Redistricting of Arizona with Population Distribution}
\end{table}

With one million simulations of our algorithm, we obtain an ensemble of one million potential redistricting solutions for North Carolina. We preserve all the county borders except for the three split counties. To reach exact population equity in a district plan, the population percentage of state population in each district would be $100\%/13\approx 7.69\%$. To see the variation across district plans in our ensemble, we compute for each district plan the standard deviation of the differences in population between the district populations and the ideal district population, with results ranging from $0.32$ percentage points to $2.33$ percentage points with a mean of $1.11$ percentage points. Then, we extract the solutions that have a relatively low population variation among the districts. We consider the solutions that have less than or equal to $1$ percentage point of state population standard deviation as ``good'' solutions. We recognize that here we are taking a different approach than that mandated currently by HB92. As discussed earlier, we point out that many confounding factors exist (e.g., age distribution, voter registration levels, voter turnout) so that even if a district map reaches the population equity requirement of HB92, it is likely that the number of votes cast varies more than the 0.1\% from ideal that HB92 stipulates. We believe that strong results for other measures can plausibly be paired with reasonable population equity levels.

Among one million solutions obtained from our simulations, there is a total of $533$ duplicate results, which is about $0.05\%$ of the total number of solutions. It is interesting to see that the number of duplicates is quite low. After removing the duplicate solutions, we have a total of $330,493$ ``good'' solutions. 

In Figures~\ref{fig:Sample Result 1 on Potential Redistricting plan for North Carolina from MATLAB} and~\ref{fig:Sample Result 2 on Potential Redistricting plan for North Carolina from MATLAB}, two samples of redistricting solutions in the ensemble supplied by the algorithm are shown. Both $part1(:,:,1)$ and $part1(:,:,2)$ are proposed redistricting plans in which each entry (except $0$) represents a generalized county and each column represents a district. For display purposes, we have truncated these matrices not to show rows consisting of all zeros. The best solution (in terms of smallest population standard deviation, $0.32$ percentage points) from the simulations is shown in Table~\ref{tab: Best Solution of the Simulation for Redistricting of Arizona with Population Distribution}.

\section{Analysis}\label{sec:analysis}

In this section, we apply several gerrymandering metrics to our ensemble of district plans and analyze the results.

\subsection{Metrics Used}

\subsubsection{Compactness}

Compactness serves as a crucial criterion in redistricting, often used to identify partisan gerrymandering. In North Carolina's HB92, perimeter compactness is highlighted, aligning with our choice of the Polsby-Popper measure. Notably, both the Polsby-Popper and Reock scores were considered in the North Carolina case \textit{Common Cause v. Rucho} (2018), emphasizing their relevance in evaluating district plans.

Denote $\xi = \{D_1, \ldots, D_{13}\}$ as a district plan containing information for the $13$ districts.  Compactness metrics are defined per district, and are consequently used to compute both the minimum and average of metric over all $D_i \in \xi$.  Minimum corresponds to the worst district for that metric, and average as a measurement over the entire district plan $\xi$.  Egregious values in either case can indicate partisan gerrymandering, since it only takes one district to change the number of elected representatives, and an overall score violation implies multiple districts to be affected.




A common measure of compactness is the \emph{Polsby-Popper} measure, which is defined as
\begin{equation}
    PP(D_i) = 4 \pi\cdot\frac{\text{area of district } i}{\left(\text{perimeter of district } i\right)^2}.\label{eq:pp}
\end{equation}

For a district that is a perfect circle, $PP(D_i) = 1$.  The scaling factor of $4\pi$ normalizes $PP(D_i)$ to be in $[0,1]$, where $PP(D_i)=1$ is the most compact shape based on this measure.  
An entire plan $\xi$ is often evaluated on its average score, in the NC case given by
\begin{equation}\label{eq:pp-avg}
    PP_{\text{avg}}(\xi) = \frac{1}{13} \sum_{i=1}^{13} PP(D_i)
\end{equation}

We remark that this is comparatively equivalent to the Schwartzberg compactness measure in the sense that a district with larger Polsby-Popper score always has a smaller Schwartzberg score and vice versa.







\subsubsection{Partisan Measures}
\label{sec:additional-metrics}

The efficiency gap is a measure of the number of wasted partisan votes in an election.  Note that in a purely bipartisan election, 50\% of votes are wasted (all for the losing candidate, and anything  beyond the 50\% winning threshold for the winner).  If one party has substantially more wasted votes than another, it is potential evidence toward gerrymandering.  For a given plan $\xi$, we formally define
\begin{equation}
    \text{efficiency gap}(\xi) = \frac{[\text{dem.~wasted votes}] - [\text{rep.~wasted votes}]}{\text{votes cast statewide}},
\end{equation}
where dem.~ and rep.~wasted votes refers to the overall statewide wasted votes. In a two-candidate election, this is a number between -0.5 and 0.5, where a strictly positive value indicates gerrymandering against the Democratic Party (strictly negative against Republican).  The extreme case of $\text{efficiency gap}(\xi)\approx 0.5$ implies that all Democratic votes are wasted in each $D_i \in \xi$, either in a landslide victory for Democrats, or during a Republican victory by 1 vote.  The $-0.5$ case occurs vice versa.

The \textit{mean-median difference} is a measure of skewness and simply subtracts the median (Democratic) voter percentage over the 13 districts from the mean (Democratic) voter percentage.  This is a traditional skewness measure used in statistics, and is interpreted similarly to the efficiency gap in that (large) positive values indicate wasted Democratic votes.

Lastly, the \textit{lopsided margin} measure subtracts the average Republican winning \% out of all Republican winning counties from the average Democratic winning \% out of all Democratic winning counties.  Note that before differencing, both percentages are greater than 50\%.  A large value then suggests that the Democratics are wasting votes, from their high average and/or the Republican low average, and vice versa.

\subsection{Results}\label{sec:results}



\begin{minipage}{\textwidth}
  \begin{minipage}[b]{0.49\textwidth}
    \centering
  \includegraphics[scale=0.7]{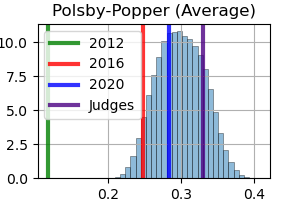}
  \captionof{figure}{Compactness Scores comparison using various measures. The histograms show results for our ensemble of district maps.}
  \label{fig:compactness}
 \end{minipage}
 \hfill
  \begin{minipage}[b]{0.49\textwidth}
    \centering
  \resizebox{\columnwidth}{!}{\begin{tabular}{|c|c|c|c|c|c|}\hline
    \multirow{1}{*}{\#} & \multicolumn{4}{c|}{Congressional District Maps} & \multirow{1}{*}{Our} \\ 
    \multirow{1}{*}{Splits} & \multicolumn{1}{c}{2012} & \multicolumn{1}{c}{2016} & \multicolumn{1}{c}{2020} & \multicolumn{1}{c|}{Judges} & \multirow{1}{*}{Alg.}\\\hline
   2 & 31 & 13 & 12 & 12 &1  \\
   3 & 8 & 0 & 0 & 0 & 0     \\
   4 & 1 & 0 & 0 & 0 &2 \\\hline \hline
   \textbf{Total} & 40 & 13 & 12 & 12 & 3 \\\hline
   \end{tabular}}
  \captionof{table}{Number of splits in Counties on North Carolina Congressional maps versus our approach. In the first column, the number of Districts into which a County might be split is indicated.  ``Our Alg.'' column refers to the maximum possible allowed splits.}
  \label{tab:Number of Splits in Counties on North Carolina Congressional Maps}
\end{minipage}
\end{minipage}

The histogram in Figure \ref{fig:compactness} provides a summary of the average Polsby-Popper score (calculated using Equation \eqref{eq:pp}) from the simulations described in Section \ref{sec:ensemble}. Overlaid on the histogram are the 2012, 2016, and 2020 North Carolina district plans, as well as the plan proposed by the judges. The results of the historic plans are clearly visible (as seen in \textit{Common Cause v. Rucho} (2018)), with the 2012 plan demonstrating clear evidence of gerrymandering on the basis of compactness.  Also consistent with historic findings is that the 2016 plan manages to produce viable compactness results despite producing a gerrymander.  Moreover, our approach shows promise for achieving compactness, with most of the results clustering around the 2020 and judges' plans, and the minimum score close to that of the 2016 plan.

Table \ref{tab:Number of Splits in Counties on North Carolina Congressional Maps} compares the number of county splits in historic and judges' plans to those generated by our algorithm.  Note that the ``Our Alg.'' column refers to the \emph{maximum allowed} splits; this informs whether a county is available for splitting or not.  This table effectively showcases one of the major advantages of our algorithm: it has only a quarter of the maximum allowable splits compared to the 12 fixed splits in the 2020 and judges' plans.

\begin{figure}
    \centering
    \includegraphics[scale=0.7]{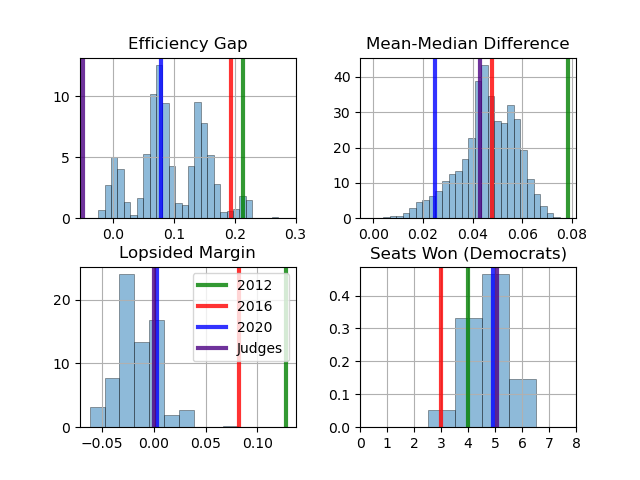}
    \caption{Partisan Scores comparison using various measures. The histograms show results for our ensemble of district maps.}
    \label{fig:partisan}
\end{figure}

In terms of partisan results, Figure \ref{fig:partisan} provides a histogram summary of the Efficiency Gap, Mean-Median Difference, and Lopsided Margin measures across simulations. To avoid gerrymandering, these measures should ideally be as close to zero as possible. As seen in the figure, historic consistency is evident, with the 2012 plan being a clear outlier in all categories (aside from Seats Won), as well as the 2016 plan (aside from Mean-Median Difference). Our simulated plans largely cluster around the 2020 plan for efficiency gap, lopsided margins, and seats won, while clustering around the judges' plan for mean-median difference and lopsided margins. However, it is worth noting that plans may pass certain criteria and fail others, as observed with the 2016 plan's performance on mean-median difference. Therefore, it is crucial to evaluate plans across multiple measures. Overall, Figure \ref{fig:partisan} provides strong evidence that most of the plans generated by our algorithm perform well in terms of each of the partisan criteria mentioned.

\section{Discussion}\label{sec:discussion}
Political gerrymandering involves a complex network of individuals, including politicians who devise plans, judges who assess political fairness, and voters who belong to communities sharing values and culture, all of whom are invested in ensuring their votes are not invalidated through gerrymandering. The role of mathematicians is also crucial, as they develop criteria for assessing gerrymandering and algorithms for constructing new plans. However, gerrymandering is fundamentally based on obfuscation and the use of esoteric algorithms.

In this case study, our algorithm, based on intuitive notions, has been shown to generally pass common gerrymandering tests and produce results similar to those of judges' and 2020 plans in North Carolina. Our algorithm has verified the gerrymanders of 2012 and 2016 in terms of their compactness and partisan measures. 

In future work on North Carolina's districting, we will still heavily minimize the number of county splits, but we will consider splitting the counties differently. Our current splitting technique is an ad hoc process that could be developed further to try to avoid potential manipulative usage of splitting of counties. For cases where additional county splits may be necessary, particularly to meet the stringent HB92 population equity mandate, county split metrics such as \emph{split pairs} introduced by \citet{waschpress2021split} could offer valuable insights. This metric could provide a more nuanced understanding of the implications of different county-splitting strategies.

Also, county demographics can be taken into consideration to ensure the contiguity of potential communities of interest. Since the ``big rocks," which are pre-selected counties, were pre-determined in our model, a few counties will not have a chance to be included in the same district. To better improve our model, we would like to use different sets of pre-selected counties and run the corresponding simulations. Although the logic and steps will be similar in each algorithm based on our model, the characteristics of the pre-selected counties and the algorithm will need to be updated. These characteristics can include county population, the population of neighboring counties, the location and the distance between pre-selected counties, and additional characteristics.

The ``Rocks-Pebbles-Sand" approach developed in this case study has useful features for usage in other states as well. Most importantly, our method is easily explainable to politicians, judges, and the general voting public, and emphasizes the importance of preserving communities of interest by minimizing splits and keeping counties intact. The ``Rocks-Pebbles-Sand" approach naturally maintains community aspects such as minority percentages, while also honoring the ``One Person-One Vote Principle" by choosing population equity to minimize at each step of the algorithm while prioritizing the maps that are relatively compact. As indicated in Appendix \ref{app:addressing-HB92} (Figure \ref{fig:population_deviation}), our method could meet a smaller population equity constraint if it was hard coded into the algorithm (e.g., 1\%, similar to the approach used in \citet{cirincione2000assessing}). Although additional splits could feasibly allow for meeting North Carolina's 0.1\% constraint, the exact number required and strategy of splitting is left for future work. Our method also provides results that can be easily fine-tuned to satisfy stricter bounds, such as shifting precincts to satisfy population equity requirements, and is easily adaptable to other states by predetermining natural county splits and building an adjacency matrix. Moreover, additional constraints based on another state legislature can be easily introduced as filtering steps in the algorithm. In summary, this case study indicates the potential of an approach built upon an easy-to-understand intuition.

\section{Conflict of Interest}
On behalf of all authors, the corresponding author states that there is no conflict of interest.

\section{Data Availability}
The data for Figures \ref{fig:historic_plans} and \ref{fig:compactness} are available at North Carolina General Assembly, \url{https://www.ncleg.gov/Redistricting} and \url{https://docs.google.com/spreadsheets/d/1XbUXnI9OyfAuhP5P3vWtMuGc5UJlrhXbzZo3AwMuHtk/edit#gid=0}. 

\vspace{3mm}
\noindent The data for Figure \ref{fig:Pre-selected Counites on Population Density map} is available at the United States Census Bureau,\url{https://www2.census.gov/geo/pdfs/reference/guidestloc/37_NorthCarolina.pdf}.
\vspace{3mm}

\noindent The data for Figures \ref{fig: 2016 U.S. House Representative Election Results of Mecklenburg County}-\ref{fig: 2016 U.S. House Representative Election Results of Wake County} and Tables \ref{tab: Results of the United States House of Representatives Elections in North Carolina from 2016 - 2020}-\ref{tab: Votes Redistribution of Wake County} are available at State of North Carolina Election Data, \url{https://er.ncsbe.gov/}.

\vspace{3mm}
\noindent 
The data for Table \ref{tab: Population of 107 Generalized Counties} is available at the United States Census \url{https://www.census.gov/data.html}.

\section{Code Availability}

Codes and Data sets generated for this study are available in the OneDrive folder as shown below:
\newline \url{https://livecsupomona-my.sharepoint.com/:f:/g/personal/yuz2_cpp_edu/El5BJcxWIylKhf5dO3cQZR4BsOmW4mE4v_F5uC3zwvJMeA?e=pJxUjE}

\bibliography{biblio.bib}

\appendix

\section{Addressing HB92}\label{app:addressing-HB92}

Among several criteria used to assess potential gerrymandering, \emph{population equity} quantifies the deviation of the population within a district compared to the \emph{ideal population} which is the total state population divided by the number of districts.  In particular, the population of district $D_i$ under census $C$ is defined as
\begin{equation}
    PE_{C}(D_i) = \frac{\left\|\textrm{Pop}_{C}(D_i) - \textrm{Pop}_{C,\textrm{ideal}}\right\|}{\text{Pop}_{C,\text{ideal}}},
\end{equation}
where $\text{Pop}_C(D_i)$ is the number of registered voters in $D_i$ as accounted for by census $C$ (e.g., $C=2010, 2020, ...$), and $\text{Pop}_{C,\text{ideal}} = \frac{1}{13} \sum_{i=1}^{13} \text{Pop}_{C}(D_i)$.  Note that in our application, $C = 2010$.  \citet{herschlag2020quantifying} utilizes a \emph{root mean squared population deviation measure}

\begin{equation}\label{eq:RMSPD}
    \text{RMSPD}_C(\xi) = \sqrt{\frac{1}{13} \sum_{i=1}^{13} \left(PE_{\text{C}}(D_i)\right)^2},
\end{equation}

\begin{figure}[h!]
    \centering
    \includegraphics[width=0.8\textwidth]{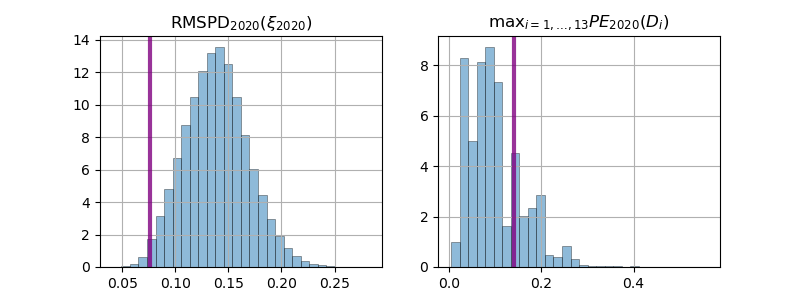}
    \caption{Root mean squared population deviation and maximum population deviation for our algorithm results. The purple vertical line indicates the 2020 plan (based on the 2010 census) when using the up-to-date 2020 census data.}
    \label{fig:population_deviation}
\end{figure}

\noindent and HB92 requires $\max_{i=1, \ldots, 13} PE_{\text{C}}(D_i) < 0.001$, that is, the maximal amount any district can deviate from its ideal population is 0.1\%.  On one hand, the HB92 principle is just in adherence with the ``One Person--One Vote Principle."  On the other, confounding factors exist within census population values (e.g., voter registration levels, voter turnout) so that the number of votes cast in a VTD is not fully explained by population representation within the VTD.  This is especially true in years toward the end of a decade.  Note that the official North Carolina 2020 plan used on the $C=2010$ census shifts precinct boundaries to produce \emph{perfect} population deviation $\text{Pop}_{2010}(D_i) = 0$ for all $i$ (up to rounding).  However, when the same plan is applied to the census gathered in the year it was used ($C=2020$), one finds $\text{RMSPD}_{2020}(\xi_{2020}) = 0.07$ and $\max_{i=1, \ldots, 13} PE_{\text{C}}(D_i) = 0.14$, 140 times larger than the requirement set by HB92! Figure \ref{fig:population_deviation} displays histograms of $\text{RMSPD}_{2020}(\xi)$ and $\max_{i=1, \ldots, 13} PE_{\text{C}}(D_i)$ for the district plans generated by our algorithm, with the $C=2020$ values overlayed for $\xi_{2020}$.


\end{document}